\documentclass[reqno,10pt]{amsart}
\usepackage{amsmath,amssymb,amsthm,graphicx,multirow,amsxtra}

\newtheorem{theorem}{Theorem}[section]

\newtheorem{corollary}[theorem]{Corollary}
\newtheorem{lemma}[theorem]{Lemma}
\newtheorem{prop}[theorem]{Proposition}
\theoremstyle{definition}
\newtheorem{definition}[theorem]{Definition}

\newtheorem{remark}[theorem]{Remark}
\numberwithin{equation}{section}

\def\({\left(}              \def\){\right)}
\def\[{\left[}              \def\]{\right]}
\newcommand{\R}{\mathbb{R}}
\newcommand{\Z}{\mathbb{Z}}
\renewcommand{\H}{\mathbb{H}}
\newcommand{\T}{\mathbb{T}}     \newcommand{\C}{\mathbb{C}}
\newcommand{\N}{\mathbb{N}}     \def\I{{\mathbb I}}
\newcommand{\D}{\Delta}
\newcommand{\s}{\sigma}
\newcommand{\e}{\varepsilon}         \newcommand{\m}{\mu}
\def\t{\tau}                \def\i{\overset{\circ}{\imath}}
  \def\b{\beta}               \def\a{\alpha}
  \def\g{\gamma}  \def\l{\lambda}

\def\P{{\mathbb P}}        
        
    \def\Log{\textup{Log}}
      \def\diag{\textup{diag}}

\def\mumax{\m_{\textup{max}}}
\begin{document}

\title[A Unified Floquet Theory]{A Unified Floquet
Theory for Discrete, Continuous, and Hybrid Periodic Linear
Systems}
\author[DaCunha]{Jeffrey J. DaCunha}%
\address{Lufkin Automation, Midland, Texas 79703 USA}
\email{jeffrey\_dacunha@yahoo.com}
\author[Davis]{John M. Davis*}%
\address{Baylor University, Department of Mathematics, Waco, Texas 76798 USA}
\email{John\_M\_Davis@baylor.edu}
\subjclass{34K13, 34C25, 39A11, 39A12}%
\keywords{Floquet theory, Floquet decomposition,  time scale,
hybrid system, Lyapunov transformation, monodromy operator, Floquet
multiplier, Floquet
exponent, stability}%
\thanks{*Work supported by grants NSF EHS\#0410685 and NSF CMMI\#0726996.}


\begin{abstract}
In this paper, we study periodic linear systems on periodic time
scales which include not only discrete and continuous dynamical
systems but also systems with a mixture of discrete and
continuous parts (e.g.~hybrid dynamical systems). We develop a
comprehensive Floquet theory including Lyapunov transformations
and their various stability preserving properties, a unified
Floquet theorem which establishes a canonical Floquet
decomposition on time scales in terms of the generalized
exponential function, and use these results to study homogeneous
as well as nonhomogeneous periodic problems. Furthermore, we
explore the connection between Floquet multipliers and Floquet
exponents via monodromy operators in this general setting and
establish a spectral mapping theorem on time scales. Finally, we
show this unified Floquet theory has the desirable property that
stability characteristics of the original system can be
determined via placement of an associated (but time varying)
system's poles in the complex plane. We include several examples
to show the utility of this theory.
\end{abstract}

\maketitle

\tableofcontents

\section{Introduction}

It is widely known that the stability characteristics of a
nonautonomous $p$-periodic linear system of differential or
difference equations can be characterized completely by a
corresponding autonomous linear system of differential or
difference equations by a periodic Lyapunov transformation of
variables \cite{Chicone, KePe, Rugh}. One application of Lyapunov
transformations has been in generating different state variable
descriptions of linear time invariant systems because different
state variable descriptions correspond to different and perhaps
more advantageous points of view in determining the system's
output characteristics. This is useful in signals and systems
applications for the simple fact that different descriptions of
state variables allow usage of linear algebra to design and study
the internal structure of a system. Having the ability to change
the internal structure without changing the input-output behavior
of the system is useful for identifying implementations of these
systems that optimize some performance criteria that may not be
directly related to input-output behavior, such as numerical
effects of round-off error in a computer-based systems
implementation.  For example, using a transformation of variables
on a discrete time non-diagonal $2\times 2$ system, one can obtain
a diagonal system matrix which separates the state update into two
decoupled first-order difference equations, and, because of its
simple structure, this form of the state variable description is
very useful for analyzing the system's properties \cite{SS}.

Without question, the study of periodic systems in general and
Floquet theory in particular has been central to the differential
equations theorist for some time. Researchers have explored these
topics for ordinary differential equations \cite{Chicone, Fr,
GeWe, Jo, Pan, PaKr, SiWuPe, Weikard}, partial differential
equations \cite{CaOtRo,ChLuMP,GeWe,Ku}, differential-algebraic
equations \cite{Demir, LaMaWi}, and discrete dynamical systems
\cite{Ahlbrandt,KePe,Tep}. Certainly \cite{MPSe} is a landmark
paper in the area. Not surprisingly, Floquet theory has wide
ranging effects, including extensions from time varying linear
systems to time varying nonlinear systems of differential
equations of the form $x'=f(t,x)$, where $f(t,x)$ is smooth and
$\omega$-periodic in $t$.  The paper by Shi \cite{Shi} ensures the
global existence of solutions and proves that this system is
topologically equivalent to an autonomous system $y'=g(y)$ via an
$\omega$-periodic transformation of variables. The theory has also
been extended by Weikard \cite{Weikard} to nonautonomous linear
systems of the form $\dot{z}=A(x)z$ where $A:\C\to\C^{n\times n}$
is an $\omega$-periodic function in the complex variable $x$,
whose solutions are meromorphic.  With the assumption that $A(x)$
is bounded at the ends of the period strip, it is shown that there
exists a fundamental solution of the form $P(x)e^{Jx}$ with a
certain constant matrix $J$ and function $P$ which is rational in
the variable $e^{2\pi i x/\omega}$.

In a relatively recent paper by Teplinski\u{i} and Teplinski\u{i}
\cite{Tep}, Lyapunov transformations and discrete Floquet theory
are extended to countable systems in $l_\infty (\mathbb{N},\R)$.
It is proved that the countable time varying system can be
represented by a countable time invariant system provided its
finite-dimensional approximations can also be represented by time
invariant systems.

Lyapunov transformations and Floquet theory have also been used to
analyze the stability characteristics of quasilinear systems with
periodically varying parameters.  In 1994, Pandiyan and Sinha
\cite{Pan} introduced a new technique for the investigation of
these systems based on the fact that all quasilinear periodic
systems can be replaced by similar systems whose linear parts are
time-invariant, via the well known Lyapunov-Floquet
transformation.

In the paper by Demir \cite{Demir}, the equivalent of Floquet
theory is developed for periodically time-varying systems of
linear DAEs: $\frac{d}{dt}(C(t)x)+G(t)x=0$ where the $n\times n$
matrices $C(\cdot)$ (not full rank in general) and $G(\cdot)$ are
periodic.  This result is developed for a direct application to
oscillators which are ubiquitous in physical systems:
gravitational, mechanical, biological, etc., and especially in
electronic and optical ones. For example, in radio frequency
communication systems, they are used for frequency translation of
information signals and for channel selection. Oscillators are
also present in digital electronic systems which require a time
reference, i.e., a clock signal, in order to synchronize
operations.  All physical systems, and in particular electronic
ones, are corrupted by undesired perturbations such as random
thermal noise, substrate and supply noise, etc. Hence, signals
generated by practical oscillators are not perfectly periodic.
This performance limiting factor in electronic systems is also
analyzed in \cite{Demir} and a theory for nonlinear perturbation
analysis of oscillators described by a system of DAEs is
developed.

The intent of this paper is to extend the current results of
continuous and discrete Floquet theory to the more general case of
an arbitrary periodic \textit{time scale} (a special case of a {\it
measure chain} \cite{Hi1, Hi2}), defined as any closed subset of
$\R$. In particular, the main result shows that there exists a not
necessarily constant $n\times n$ matrix $R(t)$ such that
$e_R(t_0+p,t_0)=\Phi_A(t_0+p,t_0)$, where $\Phi_A(t,t_0)$ is the
transition matrix for the $p$-periodic system $x^\D(t)=A(t)x(t)$,
and that the transition matrix can be represented by the product of
a $p$-periodic Lyapunov transformation matrix and a time scale
matrix exponential, i.e. $\Phi_A(t,t_0)=L(t)e_R(t,t_0)$, which is
known as the {\it Floquet decomposition} of the transition matrix
$\Phi_A(t,t_0)$.  Even though the matrix $R(t)$ is in general time
varying, because of its construction, one can still analyze the
stability of the original $p$-periodic system by the eigenvalues of
$R(t)$, just as in the familiar (and special) cases of discrete and
continuous Floquet theory. We make a note that with this unified
Floquet theory, the matrix $R(t)$ does become a constant matrix $R$
in the cases that the time scale is $\R$ or $\Z$, as one would
expect.

There has been work on generalizing the Floquet decomposition
to the time scales case in the very nice paper by Ahlbrandt and Ridenhour \cite{Ahlbrandt}.
However, there are some important distinctions between their work
and this one.  First, Ahlbrandt and Ridenhour use a different
definition of a periodic time scale. Furthermore---and very
importantly---their Floquet decomposition theorem employs the usual
exponential function whereas our approach is in terms of the generalized time
scale exponential function.  Finally, we go on to develop a complete
Floquet theory including Lyapunov transformations and stability,
Floquet multipliers, Floquet exponents, and apply this theory to questions of stability of periodic linear systems.

Additionally, this paper develops a solution to the previously open
question of the existence of a solution matrix $R(t)$ to the matrix
equation $e_R(t,\t)=M$, where $R(t)$ and $M$ are $n\times n$
matrices, and $M$ is constant and nonsingular.  Whereas the question
of a generalized time scales logarithm remains open, the solution
contained within this paper offers a step towards fulfilling
this gap in time scales analysis.

This paper is organized as follows.  In Section~\ref{LyapTrans} the
generalized Lyapunov transformation for time scales is developed and
it is shown that the change of variables using the time scales
version of this transformation preserves the stability properties of
the system. Then in Section~\ref{FloquetTheory}, the notion of a
periodic time scale is presented and the main theorem, the unified
and extended version of Floquet's theorem, is introduced for the
homogeneous and nonhomogeneous cases of a periodic system on a
periodic time scale.  Three examples are given in
Section~\ref{Examples} to show that the unified theory of Floquet is
functional in the cases $\T=\R$, $\T=\Z$, and more interestingly,
when $\T=\P_{1,1}$. Section~\ref{FloquetMultipliersSection}
introduces the unified theorems for Floquet multipliers, Floquet
exponents, as well as a generalized spectral mapping theorem for
time scales.  In Section~\ref{ExamplesRevisited}, the examples from
Section~\ref{Examples} are revisited and the theorems introduced in
Section~\ref{FloquetMultipliersSection} are illustrated.  We end the
paper with the Conclusions, where the main ideas are stated and the
results of the paper are summarized. In the Appendix, the theory of
time scales is introduced and necessary definitions and results are
stated to keep the paper relatively self-contained. An excellent
introduction to the subject of time scales analysis can be found in
Bohner and Peterson's introductory texts \cite{BoPe1,BoPe}.

\section{The Lyapunov Transformation and
Stability}\label{LyapTrans}

We begin by analyzing the stability preserving property associated
with a change of variables using a Lyapunov transformation on the
regressive time varying linear dynamic system
    \begin{equation}\label{A}
    x^\D(t)=A(t)x(t),\quad x(t_0)=x_0.
    \end{equation}

\begin{definition}
A \textit{Lyapunov transformation} is an invertible matrix
$L(t)\in \textup{C}_{\textup{rd}}^1(\T,\R^{n\times n})$ with the
property that, for some positive $\eta,\rho\in\R$,
    \begin{equation}\label{Lyaptransbound}
    ||L(t)||\leq\rho\hskip.25in\text{and}\hskip.25in \det L(t)\geq\eta
    \end{equation}
for all $t\in\T$.
\end{definition}

The two following lemmas can be found in the classic text by A.C.
Aitken \cite{Aitken}.

\begin{lemma}\label{matrixinvbd}
Suppose that $A(t)$ is an $n\times n$ matrix such that $A^{-1}(t)$
exists for all $t\in\T$.  If there exists a constant $\a>0$ such
that $||A^{-1}(t)||\leq\a$ for each $t$, then there exists a
constant $\b$ such that $|\det A(t)|\geq\b$ for all $t\in\T$.
\end{lemma}

\begin{lemma}\label{matrixbd}
Suppose that $A(t)$ is an $n\times n$ matrix such that $A^{-1}(t)$
exists for all $t\in\T$.  Then
   $$
    ||A^{-1}(t)||\leq\frac{||A(t)||^{n-1}}{|\det A(t)|}
   $$
for all $t\in\T$.
\end{lemma}

A consequence of Lemma~\ref{matrixinvbd} and Lemma~\ref{matrixbd}
is that the inverse of a Lyapunov transformation is also bounded.
An equivalent condition to \eqref{Lyaptransbound} is that there
exists a $\rho>0$ such that
    $$
    ||L(t)||\leq\rho\hskip.25in\text{and}\hskip.25in
    ||L^{-1}(t)||\leq\rho
    $$
for all $t\in\T$.
\begin{definition}
The time varying linear dynamic equation \eqref{A} is called
\textit{uniformly stable} if there exists a finite positive
constant $\gamma$ such that for any $t_0,\:x(t_0)$ the
corresponding solution satisfies
   $$
    ||x(t)||\leq\gamma ||x(t_0)||, \qquad t\geq t_0.
   $$
\end{definition}

Uniform stability can also be characterized using the following
theorem.

\begin{theorem}\label{USchar}
The time varying linear dynamic equation \eqref{A} is uniformly
stable if and only if there exists a $\gamma>0$ such that the
transition matrix $\Phi_A$ satisfies $||\Phi_A(t,t_0)||\leq\gamma$
for all $t\geq t_0$ with $t,\:t_0\in\T$.
\end{theorem}

\begin{definition} The time varying linear dynamic
equation \eqref{A} is called \textit{uniformly exponentially
stable} if there exists finite positive constants $\gamma,\:\l$
with $-\l\in\mathcal{R}^+$ such that for any $t_0,\:x(t_0)$ the
corresponding solution satisfies
    $$
    ||x(t)||\leq||x(t_0)||\gamma e_{-\l}(t,t_0),\qquad t\geq t_0.
    $$
\end{definition}

Uniform exponential stability can also be characterized using the
following theorem.

\begin{theorem}\label{UESchar}
The time varying linear dynamic equation \eqref{A} is uniformly
exponentially stable if and only if there exists an
$\l,\:\gamma>0$ with $-\l\in\mathcal{R}^+$ such that the
transition matrix $\Phi_A$ satisfies
    $$
    ||\Phi_A(t,t_0)||\leq\gamma e_{-\l}(t,t_0)
    $$
for all $t\geq t_0$ with $t,\:t_0\in\T$.
\end{theorem}

The last stability definition given uses a uniformity condition to
conclude exponential stability.
\begin{definition}\label{asymdef}
The linear state equation \eqref{A} is defined to be
uniformly asymptotically stable if it is uniformly stable
and given any $\delta > 0$, there exists a $T>0$ so that for any
$t_0\text{ and }x(t_0)$, the corresponding solution $x(t)$ satisfies
\begin{equation}\label{asym}
||x(t)||\leq \delta ||x(t_0)||,\quad t\geq t_0+T.
\end{equation}
\end{definition}
It is noted that the time $T$ that must pass before the norm of
the solution satisfies \eqref{asym} and the constant $\delta>0$ is
independent of the initial time $t_0$.

\begin{theorem}\label{Thm-UASiffUES}
The linear state equation \eqref{A} is uniformly exponentially stable if and only if
it is uniformly asymptotically stable.
\end{theorem}
\begin{proof}
Suppose that the system \eqref{A} is uniformly exponentially stable.  This implies
that there exist constants $\gamma,\lambda>0$ with $-\lambda\in\mathcal{R}^+$ so that $||\Phi_A(t,\t)||\leq\gamma
e_{-\lambda}(t,\t)$ for $t\geq\t$. Clearly, this implies uniform stability.  Now, given a $\delta>0$, we choose a sufficiently
large positive constant $T>0$ such that $t_0+T\in\T$ and $e_{-\lambda}(t_0+T,t_0)\leq\frac{\delta}{\gamma}$. Then for any $t_0$
and $x_0$, and $t\geq T+t_0$ with $t,T+t_0\in\T$,
\begin{align*}
||x(t)||&=||\Phi_A(t,t_0)x_0||\\
&\leq ||\Phi_A(t,t_0)||\:||x_0||\\
&\leq\gamma e_{-\lambda}(t,t_0)||x_0||\\
&\leq\gamma e_{-\lambda}(t_0+T,t_0)||x_0||\\
&\leq\delta||x_0||,\quad t\geq t_0+T.
\end{align*}Thus, \eqref{A} is uniformly asymptotically stable.

Now suppose the converse.  By definition of uniform asymptotic stability, \eqref{A}
is uniformly stable.  Thus, there exists a constant $\gamma>0$ so that
\begin{equation}\label{16}
||\Phi_A(t,\t)||\leq\gamma,\quad \text{ for all } t\geq\t.
\end{equation}
Choosing $\delta = \frac{1}{2}$, let $T$ be a positive constant so that $t_0+T\in\T$ and \eqref{asym} is satisfied.  Given a
$t_0$ and letting $x_a$ be so that $||x_a||=1$, we have
$$
||\Phi_A(t_0+T,t_0)x_a||=||\Phi_A(t_0+T,t_0)||.
$$
When $x_0=x_a$, the solution $x(t)$ of \eqref{A} satisfies
$$
||x(t_0+T)||=||\Phi_A(t_0+T,t_0)x_a||=||\Phi_A(t_0+T,t_0)||\:||x_a||\leq\frac{1}{2}||x_a||.
$$
From this, we obtain
\begin{equation}\label{17}
||\Phi_A(t_0+T,t_0)||\leq\frac{1}{2}.
\end{equation}
It can be seen that for any $t_0$ there exists an $x_a$ as
claimed. Therefore, the above inequality~\eqref{17} holds for any $t_0\in\T$.
%
%
%
%
%

Using the bound from the inequalities \eqref{16} and \eqref{17}, we have the following set of inequalities on intervals in
the time scale of the form $\[\t+kT,\t+(k+1)T\)_\T$, with arbitrary $\t$:
\begin{align*}
||\Phi_A(t,\t)||&\leq \g,\quad t\in\[\t,\t+T\)_\T\\
||\Phi_A(t,\t)||&=||\Phi_A(t,\t+T)\Phi_A(\t+T,\t)||\\
&\leq||\Phi_A(t,\t+T)||\:||\Phi_A(\t+T,\t)||\\
&\leq\frac{\g}{2},\quad t\in\[\t+T,\t+2T\)_\T\\
||\Phi_A(t,\t)||&=||\Phi_A(t,\t+2T)\Phi_A(\t+2T,\t+T)\Phi_A(\t+T,\t)||\\
&\leq||\Phi_A(t,\t+2T)||\:||\Phi_A(\t+2T,\t+T)||\:||\Phi_A(\t+T,\t)||\\
&\leq\frac{\g}{2^2},\quad t\in\[\t+2T,\t+3T\)_\T.\\
\end{align*}
In general, for any $\t\in\T$, we have
$$
||\Phi_A(t,\t)||\leq\frac{\g}{2^k},\quad
t\in\[\t+kT,\t+(k+1)T\)_\T.
$$
We now choose the bounds to obtain a decaying exponential bound. We now find a positive constant $\l$
(with $-\l\in\mathcal{R}^+$) that satisfies $\frac{1}{2}\leq e_{-\l}(\t+T,\t)\leq e_{-\l}(t,\t)$, for $t\in\[\t,\t+T\)_\T$.
Recall the fact 
for any positive constant $\b$ such that $-\b\in\mathcal{R}^+$,
$$
\int_{t_0}^t\lim_{s\searrow\mumax}\xi_s(-\b)\D\t\leq\int_{t_0}^t\lim_{s\searrow\mu(\t)}\xi_s(-\b)\D\t.
$$
Then we must define $\l$ such that
$$
\ln\left(\frac{1}{2}\right)=\int_{t_0}^{t_0+T}\lim_{s\searrow\mumax}\xi_s(-\l)\D\t.
$$
It is easy to verify that
$$
\l = \lim_{s\searrow\mumax}\frac{1-\frac{1}{2}^\frac{s}{T}}{s}.
$$
Since $T$ has been established to satisfy~\eqref{asym} with $\delta=\frac{1}{2}$, for any $k\in\N_0$ and $t\in\[\t+kT,\t+(k+1)T\)_\T$, where $\t+kT$ and $\t+(k+1)T$ are merely the
upper bounds on the interval, not necessarily elements in the time scale, we now have
$$
\frac{1}{2^k}\leq e_{-\l}(t,\t),
$$
for $t\in\[\t+kT,\t+(k+1)T\)_\T$.

Then for all $t,\t\in\T$ with $t\geq\t$, we obtain the correct value of $\l$ for the decaying exponential bound
$$
||\Phi_A(t,\t)||\leq\gamma e_{-\lambda}(t,\t).
$$
Therefore, by Theorem~\ref{UESchar} we have uniform exponential
stability.
\end{proof}
\begin{theorem}\label{GG}
Suppose $L(t)\in \textup{C}_{\textup{rd}}^1(\T,\R^{n\times
n})$, with $L(t)$ invertible for all $t\in\T$ and $A(t)$ is from
the dynamic linear system \eqref{A}. Then the transition matrix
for the system
    \begin{equation}\label{Gsys}
    Z^\D(t)=G(t)Z(t),\quad Z(\t)=I
    \end{equation}
where
    \begin{equation}\label{G}
    G(t)=L^{\s^{-1}}(t)A(t)L(t)-L^{\s^{-1}}(t)L^\D(t)
    \end{equation}
is given by
    $$
    \Phi_G(t,\t)=L^{-1}(t)\Phi_A(t,\t)L(\t).
    $$
for any $t,\t\in\T$.
\end{theorem}

\begin{proof}
First we see that by definition, $G(t)\in
\textup{C}_{\textup{rd}}(\T,\R^{n\times n})$.  For any $\t\in\T$,
we define
    \begin{equation}\label{X}
    X(t)=L^{-1}(t)\Phi_A(t,\t)L(\t).
    \end{equation}
It is obvious that for $t=\t$, $X(\t)=I$. Temporarily rearranging
\eqref{X} and differentiating $L(t)X(t)$ with respect to $t$, we
obtain \cite[Theorem 5.3 (iv)]{BoPe}
    $$
    L^\D(t)X(t) + L^\s(t)X^\D(t)=\Phi_A^\D(t,\t)L(\t)=A(t)\Phi_A(t,\t)L(\t),\\
    $$
and thus
    \begin{align*}
    L^\s(t)X^\D(t)&=A(t)\Phi_A(t,\t)L(\t)-L^\D(t)X(t)\\
    &=A(t)\Phi_A(t,\t)L(\t)-L^\D(t)L^{-1}(t)\Phi_A(t,\t)L(\t)\\
    &=[A(t)-L^\D(t)L^{-1}(t)]\Phi_A(t,\t)L(\t).
    \end{align*}
Multiplying both sides by $L^{\s^{-1}}(t)$ and noting \eqref{G}
and \eqref{X},
    \begin{align*}
    X^\D(t)&=[L^{\s^{-1}}(t)A(t)-L^{\s^{-1}}(t)L^\D(t)L^{-1}(t)]\Phi_A(t,\t)L(\t
    )\\
    &=[L^{\s^{-1}}(t)A(t)L(t)-L^{\s^{-1}}(t)L^\D(t)]L^{-1}(t)\Phi_A(t,\t)L(\t)\\
    &=G(t)X(t).
    \end{align*}
This is valid for any $\t\in\T$.  Thus, the transition matrix of
$X^\D(t)=G(t)X(t)$ is $\Phi_G(t,\t)=L^{-1}(t)\Phi_A(t,\t)L(\t)$.
Additionally, if the initial value specified in \eqref{Gsys} was
not the identity, i.e. $Z(t_0)=Z_0\neq I$, then the solution is
$X(t)=\Phi_G(t,\t)Z_0$.
\end{proof}

\subsection{Preservation of Uniform Stability}

\begin{theorem}
Suppose that $z(t)=L^{-1}(t)x(t)$ is a Lyapunov transformation.
Then the system \eqref{A} is uniformly stable if and only if
    \begin{equation}\label{L1}
    z^\D(t)=\left[L^{\s^{-1}}(t)A(t)L(t)-L^{\s^{-1}}(t)L^\D(t)\right]z(t),\quad
    z(t_0)=z_0
    \end{equation}
is uniformly stable.
\end{theorem}

\begin{proof}
Equation \eqref{A} and equation \eqref{L1} are related by the
change of variables $z(t)=L^{-1}(t)x(t)$.  By Theorem~\ref{GG},
the relationship between the two transition matrices is
    $$
    \Phi_G(t,t_0)=L^{-1}(t)\Phi_A(t,t_0)L(t_0).
    $$

Suppose that \eqref{A} is uniformly stable.  Then there exists a
$\gamma> 0$ such that $||\Phi_A(t,t_0)||\leq\gamma$ for all
$t,t_0\in\T$ with $t\geq t_0$.  Then by Lemma~\ref{matrixbd} and
Theorem~\ref{USchar}, we have
    \begin{align*}
    ||\Phi_G(t,t_0)||&=||L^{-1}(t)\Phi_A(t,t_0)L(t_0)||\\
    &\leq ||L^{-1}(t)||\:||\Phi_A(t,t_0)||\:||L(t_0)||\\
    &\leq\frac{\gamma\rho^n}{\eta}=\gamma_G,
    \end{align*}
for all $t,t_0\in\T$ with $t\geq t_0$.  By Theorem~\ref{USchar},
since $||\Phi_G(t,t_0)||\leq\gamma_G$, the system \eqref{L1} is
uniformly stable. The converse is similar.
\end{proof}

\subsection{Preservation of Uniform Exponential Stability}

\begin{theorem}
Suppose that $z(t)=L^{-1}(t)x(t)$ is a Lyapunov transformation.
Then the system \eqref{A} is uniformly exponentially stable if and
only if
    \begin{equation}\label{L2}
    z^\D(t)=\left[L^{\s^{-1}}(t)A(t)L(t)-L^{\s^{-1}}(t)L^\D(t)\right]z(t),\quad
    z(t_0)=z_0
    \end{equation}
is uniformly exponentially stable.
\end{theorem}

\begin{proof}
Equations \eqref{A} and \eqref{L2} are related by the change of
variables $z(t)=L^{-1}(t)x(t)$.  By Theorem~\ref{GG}, the
relationship between the two transition matrices is
    $$
    \Phi_G(t,t_0)=L^{-1}(t)\Phi_A(t,t_0)L(t_0).
    $$

Suppose that \eqref{A} is uniformly exponentially stable. Then
there exists an $\l,\:\gamma>0$ with $-\l\in\mathcal{R}^+$ such
that $||\Phi_A(t,t_0)||\leq\gamma e_{-\l}(t,t_0)$ for all $t\geq
t_0$ with $t,t_0\in\T$.  Then by Lemma~\ref{matrixbd} and
Theorem~\ref{UESchar}, we have
    \begin{align*}
    ||\Phi_G(t,t_0)||&=||L^{-1}(t)\Phi_A(t,t_0)L(t_0)||\\
    &\leq ||L^{-1}(t)||\:||\Phi_A(t,t_0)||\:||L(t_0)||\\
    &\leq\frac{\gamma\rho^n}{\eta}e_{-\l}(t,t_0)=\gamma_Ge_{-\l}(t,t_0),
    \end{align*}
for all $t,t_0\in\T$ with $t\geq t_0$.

By Theorem~\ref{UESchar}, since
$||\Phi_G(t,t_0)||\leq\gamma_Ge_{-\l}(t,t_0)$, the system
\eqref{L2} is uniformly exponentially stable. The converse is
similar.
\end{proof}

\begin{corollary}
Suppose that $z(t)=L^{-1}(t)x(t)$ is a Lyapunov transformation.
Then the system \eqref{A} is uniformly asymptotically stable if
and only if
    $$
z^\D(t)=\left[L^{\s^{-1}}(t)A(t)L(t)-L^{\s^{-1}}(t)L^\D(t)\right]z(t),
\quad z(t_0)=z_0
    $$
is uniformly asymptotically stable.
\end{corollary}

\begin{proof}
The proof follows from Theorem~\ref{Thm-UASiffUES}.
\end{proof}

\section{A Unified Floquet Theory}\label{FloquetTheory}
We want to make a vital distinction that is necessary for comprehension
of the notation that will be used in the remainder of the paper. When considering solutions to the linear dynamic system~\eqref{A}, just as in
the familiar case of $\T=\R$, we may write the transition matrix
differently, depending on properties of the system matrix.

In the most general case, the solution may always be expressed
as $x(t)=\Phi_A(t,t_0)x_0$, where $\Phi_A$ is written as~\cite[Theorem 3.2]{DaCunha2}
    \begin{align}\label{generalPBseries}
    \Phi_A(t,t_0)= I\notag &+\int_{t_0}^tA(\t_1)\D
    \t_1+\int_{t_0}^tA(\t_1)\int_{t_0}^{\t_1}A(\t_2)\D \t_2\D \t_1\\
    &+\cdots+\int_{t_0}^tA(\t_1)\int_{t_0}^{\t_1}A(\t_2)\cdots
    \int_{t_0}^{\t_{i-1}}A(\t_i)\D \t_i\cdots \D \t_1+\cdots,
    \end{align}
which generalizes the classical derivation using Picard iterates for a first order linear differential equation to any time scale.

When the system matrix commutes with its integral, i.e.
    $$
    A(t)\int_{s}^tA(\t)\D\t=\int_{s}^tA(\t)\D\t\,A(t),\quad\text{for all }s,t\in[t_0,\infty)_\T,
    $$
we may write the solution to~\eqref{A} as $x(t)=e_A(t,t_0)x_0$.  This type of system is called a {\it Lappo-Danilevski\v{i} system} by the Russian school; see \cite{Adrianova}. In this case, the representation of the matrix exponential $e_A$ is equivalent to the transition matrix $\Phi_A$ in~\eqref{generalPBseries}. However, the property that separates the transition matrix from the matrix exponential (defined as the solution to~\eqref{A}) is that
    $$
    A(t)e_A(t,t_0)=e_A(t,t_0)A(t),\text{ but }A(t)\Phi_A(t,t_0)\neq\Phi_A(t,t_0)A(t).
    $$

Note that if the system matrix is constant, we may again
represent the solution with the matrix exponential, either using
the previous representation or the infinite
series~\cite[Theorem 4.1]{DaCunha2}
    $$
    e_A(t,t_0)=\sum_{k=0}^\infty A^k h_k(t,t_0).
    $$

Just as on $\R$ and $h\Z$, one of the main properties that distinguishes $\Phi_A$ from $e_A$ in the context of a solution to a linear dynamic system of the form~\eqref{A}, is that in general $\Phi_A$ does not commute
with $A$.  However, when the solution to~\eqref{A} may be expressed as $e_A$, it is equivalent to the fact that $e_A$ and $A$ necessarily commute for all $t\in\T$.

We now state definitions that will be used throughout the paper.

\begin{definition}
Let $p\in[0,\infty)$.  Then the time scale $\T$ is {\it
$p$-periodic} if for all $t\in\T$:
    \begin{itemize}
    \item[\textup{(i)}] $t\in\T$ implies $t+p\in\T$,
    \item[\textup{(ii)}] $\mu(t)=\mu(t+p)$.
    \end{itemize}
\end{definition}

\begin{definition}
$A:\T\to\R^{n\times n}$ is
{\it $p$-periodic} if $A(t)=A(t+p)$ for all $t\in\T$.
\end{definition}

\subsection{The Homogeneous Equation}
Let $\T$ be a $p$-periodic time scale. Consider the regressive time varying linear dynamic initial value problem
    \begin{equation}\label{F}
    x^\D(t)=A(t)x(t),\quad x(t_0)=x_0,
    \end{equation}
where $A(t)$ is $p$-periodic for all $t\in\T$ and the time scale
$\T$ is also $p$-periodic. Although it is not necessary that the
period of $A(t)$ and the period of the time scale to be equal,
we will assume so for simplicity.  Furthermore, we will assume that $0\leq\mu(t)\leq p$ and that
the linear dynamic system and the time scale under discussion are
$p$-periodic, unless references to more general systems are
made.

We begin with a definition compulsory for the derivation of the matrix that will be introduced in Theorem~\ref{defofmatrixRthm}.

\begin{definition}\label{matrix M to a real power definition}
Given an $n\times n$ nonsingular matrix $M$ and any $r\in\R$, we
define the principal value of the real power of the matrix $M$ by
    $$
    M^r=\sum_{i=1}^kP_i(M)\l_i^r\left[\sum_{j=0}^{m_i-1}\frac{\Gamma(r+1)}{j!\;\Gamma(r-j+1)}\left(\frac{M-\l_iI}{\l_i}\right)^j\right],
    $$
where $m_i$ is the multiplicity of the eigenvalue $\l_i$ and the $P_i$ are the projection matrices of $M$ (See Appendix).
\end{definition}

\begin{remark}
Only the principal values of the matrix functions and the eigenvalues will be considered in the remainder of the paper.
\end{remark}

The following theorem is crucial for the development of the unified Floquet decomposition.

\begin{theorem}[{\bf Construction of the $R$ Matrix}]\label{defofmatrixRthm}
Given a nonsingular constant matrix $M$ and constant $p>0$, a
solution of the time scale matrix exponential equation
$e_{R}(t_0+p,t_0)=M$ is given by $R:\T\to\C^{n\times n}$, where
    \begin{equation}\label{R_formula}
    R(t):=\lim_{s\searrow\m(t)}\frac{M^\frac{s}{p}-I}{s}.
    \end{equation}
If $\T$ has constant graininess
on the interval $\[t_0, t_0+p\]_\T$, then $R(t)$ is constant.
\end{theorem}

\begin{proof}
Observe that $R(t)\int_{t_0}^tR(\t)\D\t=\int_{t_0}^tR(\t)\D\t R(t)$ for all $t,t_0\in\T$.  We can
represent the transition matrix solution $\Phi_R(t,t_0)$ of the
linear dynamic system
    $$
    z^\D(t)=R(t)z(t),\quad z(t_0)=z_0,
    $$
using the matrix exponential.  In other words, due to the commutativity of $R(t)$, we have
    $$
    \Phi_R(t,t_0)\equiv e_R(t,t_0).
    $$
It can be verified by direct calculation that $e_R^\D(t,t_0)=R(t)e_R(t,t_0)=e_R(t,t_0)R(t)$ for all $t\in\T$.
Using~\eqref{R_formula}, we obtain
    \begin{align*}
    e_R(t,t_0)&=M^{\frac{t-t_0}{p}}.
    \end{align*}
To prove this claim, first note that $e_R(t_0,t_0)=M^0=I$.  Delta differentiating, we have
    \begin{align*}
    e_R^\D(t,t_0)&=\lim_{s\searrow\m(t)}\frac{M^{\frac{t+\m(t)-t_0}{p}}-M^{\frac{t-t_0}{p}}}{s}\\
    &=\left(\lim_{s\searrow\m(t)}\frac{M^{\frac{\m(t)}{p}}-I}{s}\right)M^{\frac{t-t_0}{p}}\\
    &=R(t)e_R(t,t_0).\\
    \end{align*}
It follows straightforwardly that
    $$
    e_R(t_0+p,t_0)=M^{\frac{t_0+p-t_0}{p}}=M.
    $$
\end{proof}


An interesting and useful property of the matrix $R$ constructed in Theorem~\ref{defofmatrixRthm} is stated in the following corollary.

\begin{corollary}\label{cor_evals of M are the evals of R}
The matrices $R(t)$ and $M$ have identical eigenvectors.
\end{corollary}

\begin{proof}
Implementing Theorem~\ref{matrix M to a real power theorem}, for any of the $k\leq n$ distinct eigenpairs $\{\l_i,v_i\}$ of $M$ we have
    $$
    Mv_i=\l_i v_i\implies \lim_{s\searrow\m(t)}M^{\frac{s}{p}}v_i=\lim_{s\searrow\m(t)}\l_i^{\frac{s}{p}}v_i\implies R(t)v_i=\lim_{s\searrow\m(t)}\left(\frac{\l_i^{\frac{s}{p}}-1}{s}\right)v_i.
    $$
Thus, the $k\leq n$ distinct eigenpairs of $R(t)$ are $\{\g_i(t),v_i\}_{i=1}^k$, where $\g_i(t):=\lim_{s\searrow\m(t)}\frac{\l_i^{\frac{s}{p}}-1}{s}$.
\end{proof}

\begin{lemma}\label{eRshifted}
Suppose that $\T$ is a $p$-periodic time scale and
$P(t)\in\mathcal{R}(\T,\R^{n\times n})$ is also $p$-periodic. Then
the solution of the dynamic matrix initial value problem
    \begin{equation}\label{R1}
    Z^\D(t)=P(t)Z(t),\qquad Z(t_0)=Z_0,
    \end{equation}
is unique up to a period $p$ shift.  That is,
$\Phi_P(t,t_0)=\Phi_P(t+kp,t_0+kp)$, for all $t\in\T$ and
$k\in\N_0$.
\end{lemma}

\begin{proof}
By \cite{BoPe}, the unique solution to \eqref{R1} is
$Z(t)=\Phi_P(t,t_0)Z_0$. Observe
    $$
    \Phi_P^\D(t,t_0)Z_0=P(t)\Phi_P(t,t_0)Z_0
    \qquad\text{and}\qquad
    \Phi_P(t,t_0)|_{t=t_0}Z_0=\Phi_P(t_0,t_0)Z_0=Z_0.
    $$
Now we show that $\Phi_P(t,t_0)=\Phi_P(t+kp,t_0+kp)$, for all $t\in\T$ and $k\in\N_0$.  We do so by observing that $\Phi_P(t+kp,t_0+kp)Z_0$ also solves the
matrix initial value problem \eqref{R1}.  We see that
    \begin{align*}
    \Phi_P^{\D_{t+kp}}(t+kp,t_0+kp)Z_0&=P(t+kp)\:\Phi_P(t+kp,t_0+kp)\\
    &=P(t)\:\Phi_P(t+kp,t_0+kp),\quad\text{and}\\
    \Phi_P(t+kp,t_0+kp)|_{t+kp=t_0+kp}&=\Phi_P(t+kp,t_0+kp)|_{t=t_0}\\
    &=\Phi_P(t_0+kp,t_0+kp)Z_0\\
    &=Z_0.
    \end{align*}
Since the solution to \eqref{R1} is unique, we conclude
    $$
    \Phi_P(t+kp,t_0+kp)=\Phi_P(t,t_0),\quad\text{for all }t\in\T\text{
    and }k\in\N_0.
    $$
\end{proof}



%



The next theorem is the unified and extended version of Floquet decomposition for $p$-periodic time varying linear dynamic systems.
\begin{theorem}[{\bf The Unified Floquet Decomposition}]\label{Floquet2} The transition matrix for a $p$-periodic
$A(t)$ can be written in the form
    \begin{equation}\label{Fdecomp}
    \Phi_A(t,\t)=L(t)e_R(t,\t)L^{-1}(\t)\quad\text{for all }t,\t\in\T,
    \end{equation}
where $R:\T\to\C^{n\times n}$ and $L(t)\in
\textup{C}_{\textup{rd}}^1(\T,\C^{n\times n})$ are both
$p$-periodic and invertible at each $t\in\T$.  We refer to
\eqref{Fdecomp} as the {\em Floquet decomposition for $\Phi_A$.}
\end{theorem}

\begin{proof}
We define the matrix $R$ as in Theorem~\ref{defofmatrixRthm}, with $M:=\Phi_A(t_0+p,t_0)$.  Using this definition, $R$ satisfies the equation
    $$
    e_R(t_0+p,t_0)=\Phi_A(t_0+p,t_0).
    $$
Define the matrix $L(t)$ by
    \begin{equation}\label{L}
    L(t):=\Phi_A(t,t_0)e_R^{-1}(t,t_0).
    \end{equation}
By definition, $L(t)\in
\textup{C}_{\textup{rd}}^1(\T,\C^{n\times n})$ and is invertible at each $t\in\T$. Now
    \begin{equation}\label{phi_exp_eq}
    \Phi_A(t,t_0)=L(t)e_R(t,t_0),
    \end{equation}
yields
    $$
    \Phi_A(t_0,t)=e_R^{-1}(t,t_0)L^{-1}(t)=e_R(t_0,t)L^{-1}(t),
    $$
which, together with \eqref{L}, implies
    $$
    \Phi_A(t,\t)=L(t)e_R(t,\t)L^{-1}(\t),
    $$
for all $\t,t\in\T$.

We conclude by showing that $L(t)$ is $p$-periodic.  By
\eqref{L} and Lemma~\ref{eRshifted},
    \begin{align*}
    L(t+p)&=\Phi_A(t+p,t_0)e_R^{-1}(t+p,t_0)\\
    &=\Phi_A(t+p,t_0+p)\Phi_A(t_0+p,t_0)e_R(t_0,t_0+p)e_R(t_0+p,t+p)\\
    &=\Phi_A(t,t_0)\Phi_A(t_0+p,t_0)e_R^{-1}(t_0+p,t_0)e_R(t_0,t)\\
    &=\Phi_A(t,t_0)e_R^{-1}(t,t_0)\\
    &=L(t).
    \end{align*}
\end{proof}





\begin{theorem}\label{relationship_thm}
Let $\Phi_A(t,t_0)=L(t)e_R(t,t_0)$ be a Floquet decomposition
for $\Phi_A$. Then $x(t)=\Phi_A(t,t_0)x_0$ is a solution of the
$p$-periodic system \eqref{F} if and only if
$z(t)=L^{-1}(t)x(t)$ is a solution of the system
    \begin{equation}\label{z_system}
    z^\D(t)=R(t)z(t),\quad z(t_0)=x_0.
    \end{equation}
\end{theorem}

\begin{proof}
Suppose $x(t)$ is a solution to \eqref{F}.  Then $x(t)=\Phi_A(t,t_0)x_0=L(t)e_R(t,t_0)x_0.$ Setting
    $$
    z(t):=L^{-1}(t)x(t)=L^{-1}(t)L(t)e_R(t,t_0)x_0=e_R(t,t_0)x_0,
    $$
it follows from the construction of $R(t)$ that $z(t)$ is a solution of
\eqref{z_system}.

Now suppose $z(t)=L^{-1}(t)x(t)$ is a solution of the system
\eqref{z_system}.  Then
$z(t)=e_R(t,t_0)x_0$. Set $x(t):=L(t)z(t)$. Then,
    $$
    x(t)=L(t)e_R(t,t_0)x_0=\Phi_A(t,t_0)x_0,
    $$
so $x(t)$ is a solution of \eqref{F}.
\end{proof}

\begin{corollary}\label{cor stability thru Lyap trans}
The solutions of the system~\eqref{F} are uniformly stable \textup{(}respectively, uniformly exponentially stable, asymptotically stable\textup{)} if and only if the solutions of the system~\eqref{z_system} are uniformly stable \textup{(}respectively, uniformly exponentially stable, asymptotically stable\textup{)}.
\end{corollary}
\begin{proof}
The proof follows from the fact that the systems are related by
a Lyapunov change of variables and implementation of the
appropriate stability preservation theorem in
Section~\ref{LyapTrans}.
\end{proof}

\begin{theorem}\label{eigenthm}
Given any $t_0\in\T$, there exists an initial state
$x(t_0)=x_0\neq 0$ such that the solution of \eqref{F} is
$p$-periodic if and only if at least one of the eigenvalues of
$e_R(t_0+p,t_0)=\Phi_A(t_0+p,t_0)$ is $1$.
\end{theorem}

\begin{proof}
Suppose that given an initial time $t_0$ with $x(t_0)=x_0\neq 0$,
the solution $x(t)$ is $p$-periodic.  By Theorem~\ref{Floquet2},
there exists a Floquet decomposition of $x$ given by
    $$
    x(t)=\Phi_A(t,t_0)x_0=L(t)e_R(t,t_0)L^{-1}(t_0)x_0.
    $$
Furthermore,
    $$
    x(t+p)=L(t+p)e_R(t+p,t_0)L^{-1}(t_0)x_0=L(t)e_R(t+p,t_0)L^{-1}(t_0)x_0.
    $$
Since $x(t)=x(t+p)$ and $L(t)=L(t+p)$ for each $t\in\T$, we have
    $$
    e_R(t,t_0)L^{-1}(t_0)x_0=e_R(t+p,t_0)L^{-1}(t_0)x_0,
    $$
which implies
    $$
    e_R(t,t_0)L^{-1}(t_0)x_0=e_R(t+p,t_0+p)e_R(t_0+p,t_0)L^{-1}(t_0)x_0.
    $$
Since $e_R(t+p,t_0+p)=e_R(t,t_0)$,
    $$
    e_R(t,t_0)L^{-1}(t_0)x_0=e_R(t,t_0)e_R(t_0+p,t_0)L^{-1}(t_0)x_0,
    $$
and thus
    $$
    L^{-1}(t_0)x_0=e_R(t_0+p,t_0)L^{-1}(t_0)x_0.
    $$
Since $L^{-1}(t_0)x_0\neq 0$, we see that $L^{-1}(t_0)x_0$ is an
eigenvector of the matrix $e_R(t_0+p,t_0)$ corresponding to an
eigenvalue of $1$.

Now suppose 1 is an eigenvalue of $e_R(t_0+p,t_0)$ with
corresponding eigenvector $z_0$.  Then $z_0$ is real-valued and
nonzero.  For any $t_0\in\T$, $z(t)=e_R(t,t_0)z_0$ is
$p$-periodic. Since $1$ is an eigenvalue of $e_R(t_0+p,t_0)$ with
corresponding eigenvector $z_0$ and $e_R(t+p,t_0+p)=e_R(t,t_0),$
    \begin{align*}
    z(t+p)&=e_R(t+p,t_0)z_0\\
    &=e_R(t+p,t_0+p)e_R(t_0+p,t_0)z_0\\
    &=e_R(t+p,t_0+p)z_0\\
    &=e_R(t,t_0)z_0\\
    &=z(t).
    \end{align*}
Using the Floquet decomposition from Theorem~\ref{Floquet2} and
setting $x_0:=L(t_0)z_0$, we obtain the nontrivial solution of
\eqref{F}.  Then
    $$
    x(t)=\Phi_A(t,t_0)x_0=L(t)e_R(t,t_0)L^{-1}(t_0)x_0=L(t)e_R(t,t_0)z_0=L(t)z(t
    ),
    $$
which is $p$-periodic since $L(t)$ and $z(t)$ are $p$-periodic.
\end{proof}

\subsection{The Nonhomogeneous Equation}
We now consider the nonhomogeneous uniformly regressive time varying linear
dynamic initial value problem
    \begin{equation}\label{inA}
    x^\D(t)=A(t)x(t)+f(t),\quad x(t_0)=x_0,
    \end{equation}
where $A\in\mathcal{R}(\T,\R^{n\times n})$, $f\in
\textup{C}_{\text{prd}}(\T,\R^{n\times
1})\cap\mathcal{R}(\T,\R^{n\times 1}),$ and both are $p$-periodic
for all $t\in\T$.

\begin{lemma}\label{inhomlem}
A solution $x(t)$ of \eqref{inA} is $p$-periodic if and
only if $x(t_0+p)=x(t_0)$.
\end{lemma}

\begin{proof}
Suppose that $x(t)$ is $p$-periodic.  Then by definition of a
periodic function, $x(t_0+p)=x(t_0)$.

Now suppose that there exists a solution of \eqref{inA} such that
$x(t_0+p)=x(t_0)$.  Define $z(t)=x(t+p)-x(t)$.  By assumption and
construction of $z(t)$, we have $z(t_0)=0$.  Furthermore,
    \begin{align*}
    z^\D(t)&=\[A(t+p)x(t+p)+f(t+p)\]-\[A(t)x(t)+f(t)\]\\
    &=A(t)\[x(t+p)-x(t)\]\\
    &=A(t)z(t).
    \end{align*}
By uniqueness of solutions, we see that $z(t)\equiv 0$ for all
$t\in\T$. Thus, $x(t)=x(t+p)$ for all $t\in\T$.
\end{proof}

The next theorem uses Lemma~\ref{inhomlem} to develop criteria for
the existence of $p$-periodic solutions for any $p$-periodic
vector-valued function $f(t)$.

\begin{theorem}
For all $t_0\in\T$ and for all $p$-periodic $f(t)$, there exists
an initial state $x(t_0)=x_0$ such that the solution of
\eqref{inA} is $p$-periodic if and only if there does not exist a
nonzero $z(t_0)=z_0$ and $t_0\in\T$ such that the homogeneous
initial value problem
    \begin{equation}\label{zhom}
    z^\D(t)=A(t)z(t),\quad z(t_0)=z_0,
    \end{equation}
\textup{(}where $A(t)$ is $p$-periodic\textup{)} has a
$p$-periodic solution.
\end{theorem}

\begin{proof}
The solution of \eqref{inA} is given by
    $$
    x(t)=\Phi_A(t,t_0)x_0+\int^t_{t_0}\Phi_A(t,\sigma(\t))f(\t)\D\t.
    $$
From Lemma~\ref{inhomlem}, $x(t)$ is $p$-periodic if and only if
$x(t_0)=x(t_0+p)$ which is equivalent to
    \begin{equation}\label{alg}
    \[I-\Phi_A(t_0+p,t_0)\]x_0=\int^{t_0+p}_{t_0}\Phi_A(t_0+p,\sigma(\t))f(\t)\D
    \t.
    \end{equation}
By Theorem~\ref{eigenthm}, we must show that this algebraic equation
has a solution for any initial condition $x(t_0)$ and any
$p$-periodic $f(t)$ if and only if $e_R(t_0+p,t_0)$ has no
eigenvalues equal to one.

First, suppose $e_R(t_1+p,t_1)=\Phi_A(t_1+p,t_1)$, for some
$t_1\in\T$, and suppose there are no eigenvalues equal to one.
This is equivalent to
    $$
    \det\[I-\Phi_A(t_1+p,t_1)\]\neq 0.
    $$
Since $\Phi_A$ is $p$-periodic and invertible, we obtain the
equivalent statement
    \begin{align}
    0&\neq\det\[\Phi_A(t_0+p,t_1+p)\(I-\Phi_A(t_1+p,t_1)\)\Phi_A(t_1,t_0)\]\notag\\
    &=\det\[\Phi_A(t_0+p,t_1+p)\Phi_A(t_1,t_0)-\Phi_A(t_0+p,t_0)\],\label{det2}
    \end{align}
Since $\Phi_A(t_0+p,t_1+p)=\Phi_A(t_0,t_1)$, \eqref{det2} is 
equivalent to the invertibility of $[I-\Phi_A(t_0+p,t_0)]$, for any
$t_0\in\T$. Thus, for any initial time $t_0\in\T$ and for any
$p$-periodic $f(t)$, \eqref{alg} has a solution and that solution is
    $$
    x_0=\[I-\Phi_A(t_0+p,t_0)\]^{-1}\int^{t_0+p}_{t_0}\Phi_A(t_0+p,\sigma(\t))f(\t)\D
    \t.
    $$

Now suppose \eqref{alg} has a solution for every $t_0$ and
every $p$-periodic $f(t)$.  Given an arbitrary $t_0\in\T$,
corresponding to any $f_0$, we define a
regressive $p$-periodic vector valued function $f(t)\in
\textup{C}_{\text{prd}}(\T,\R^{n\times 1})$ by
    \begin{equation}\label{f}
    f(t)=\Phi_A(\s(t),t_0+p)f_0,\quad t\in\[t_0,t_0+p\)_\T,
    \end{equation}
extending this to the entire time scale $\T$ using the
periodicity.

By construction of $f(t)$,
    $$
    \int_{t_0}^{t_0+p}\Phi_A(t_0+p,\sigma(\t))f(\t)\D\t =
    \int_{t_0}^{t_0+p}f_0\D\t = pf_0.
    $$
Thus, \eqref{alg} becomes
    \begin{equation}\label{algf}
    \[I-\Phi_A(t_0+p,t_0)\]x_0=pf_0.
    \end{equation}
For any $f(t)$ that is constructed as in
\eqref{f}, and thus for any corresponding $f_0$, \eqref{algf} has a
solution for $x_0$ by assumption.  Therefore,
    $$
    \det\[I-\Phi_A(t_0+p,t_0)\]\neq 0.
    $$
Hence, $e_R(t_0+p,t_0)=\Phi_A(t_0+p,t_0)$ has no eigenvalue of
equal to 1. By Theorem~\ref{eigenthm}, \eqref{zhom} has no
periodic solution.
\end{proof}

\section{Examples}\label{Examples}

\subsection{Discrete Time Example}\label{discrete_example_1}

Consider the time scale $\T=\Z$ and the regressive (on $\Z$) time
varying matrix
    $$
    A(t)=\left[%
    \begin{array}{cc}
      -1& \frac{2+(-1)^t}{2} \\
      \frac{2+(-1)^t}{2} & -1 \\
    \end{array}%
    \right],
    $$
which have periods of $1$ and $2$, respectively.  The transition matrix for the homogeneous periodic
linear system of difference equations
    \begin{equation}\label{discrete_example_1_dynamics}
    \D X(t)=
    \left[%
    \begin{array}{cc}
      -1& \frac{2+(-1)^t}{2} \\
      \frac{2+(-1)^t}{2} & -1 \\
    \end{array}%
    \right]X(t),
    \end{equation}
is (from straightforward calculation) given by
    $$
    \Phi_{A}(t,0)=\frac{1}{2^{t+1}}
    \left[%
    \begin{array}{cc}
      (\sqrt{3})^t+(-\sqrt{3})^t & (\sqrt{3})^{t+1}+(-\sqrt{3})^{t+1} \\
      (\sqrt{3})^{t+1}+(-\sqrt{3})^{t+1} & (\sqrt{3})^t+(-\sqrt{3})^t \\
    \end{array}%
    \right].
    $$
As in Theorem~\ref{Floquet2}, consider the equation
    $$
    e_{R}(2,0)=\Phi_A(2,0)=\frac{1}{2^{3}}
    \left[%
    \begin{array}{cc}
      (\sqrt{3})^2+(-\sqrt{3})^2 & (\sqrt{3})^{3}+(-\sqrt{3})^{3} \\
      (\sqrt{3})^{3}+(-\sqrt{3})^{3} & (\sqrt{3})^2+(-\sqrt{3})^2 \\
    \end{array}%
    \right]%
    $$
which simplifies to
    $$
    e_{R}(2,0)=(I+R)^2=\left[%
    \begin{array}{cc}
      \frac{3}{4} & 0 \\
      0 & \frac{3}{4} \\
    \end{array}%
    \right].%
    $$
Then $R(t)$ in the Floquet decomposition is given by
    \begin{equation}\label{discrete_example_R}
    R(t):=\lim_{\mu(t)\searrow
    1}\frac{1}{\mu(t)}\left(\Phi_A(2,0)^{\frac{\mu(t)}{2}}-I\right)
    =\left[%
    \begin{array}{cc}
      \frac{\sqrt{3}}{2}-1 & 0 \\
      0 & \frac{\sqrt{3}}{2}-1 \\
    \end{array}%
    \right],%
    \end{equation}
which is constant (as expected) since $\m\equiv 1$ on $\Z$. Furthermore, the $2$-periodic matrix $L(t)$
is given by
    \begin{align}
    L(t)&:=\Phi_A(t,0)e_{R}^{-1}(t,0)\notag\\
    &=\Phi_A(t,0)(I+R)^{-t}\notag\\
    &=\frac{1}{2^{t+1}}
    \left[%
    \begin{array}{cc}
      (\sqrt{3})^t+(-\sqrt{3})^t & (\sqrt{3})^{t+1}+(-\sqrt{3})^{t+1} \\
      (\sqrt{3})^{t+1}+(-\sqrt{3})^{t+1} & (\sqrt{3})^t+(-\sqrt{3})^t \\
    \end{array}%
    \right]%
    \left[%
    \begin{array}{cc}
      (\frac{\sqrt{3}}{2})^{-t} & 0 \\
      0 & (\frac{\sqrt{3}}{2})^{-t} \\
    \end{array}%
    \right]\notag\\
    &=\frac{1}{2}
    \left[%
    \begin{array}{cc}
      1+(-1)^t & \sqrt{3}+(-1)^t(-\sqrt{3})\\
      \sqrt{3}+(-1)^t(-\sqrt{3}) & 1+(-1)^t \\
    \end{array}%
    \right].\label{discrete_example_L}%
    \end{align}
Using $R$ and $L(t)$ above, it is straightforward to verify that indeed~\eqref{phi_exp_eq}, and equivalently,~\eqref{Fdecomp} both hold here.

\subsection{Continuous Time Example}\label{conts_example_1}

Consider the time scale $\T=\R$ and the time varying matrix
    $$
    A(t)=\left[%
    \begin{array}{cc}
      -1& 0 \\
      \sin(t) & 0 \\
    \end{array}%
    \right]%
    $$
which has period $2\pi$.  The
transition matrix for the homogeneous periodic linear system of
differential equations
    \begin{equation}\label{conts_example_1_dynamics}
    \dot{X}(t)=
    \left[%
    \begin{array}{cc}
      -1 & 0 \\
      \sin t & 0 \\
    \end{array}%
    \right]X(t)
    \end{equation}
is given by
    $$
    \Phi_{A}(t,0)=
    \left[%
    \begin{array}{cc}
      e^{-t} & 0 \\
      \frac{1}{2}-\frac{e^{-t}(\cos t+\sin t)}{2} & 1 \\
    \end{array}%
    \right].
    $$
As in Theorem~\ref{Floquet2}, consider the equation
    $$
    e_{R}(2\pi,0)=\Phi_A(2\pi,0)=
    \left[%
    \begin{array}{cc}
    e^{-2\pi} & 0 \\
    \frac{1}{2}-\frac{e^{-2\pi}}{2} & 1 \\
    \end{array}%
    \right].%
    $$
Then $R(t)$ in the Floquet decomposition is given by
    \begin{equation}\label{conts_example_R}
    R(t):=\lim_{\mu(t)\searrow
    0}\frac{1}{\mu(t)}\left(\Phi_A(2\pi,0)^{\frac{\mu(t)}{2\pi}}-I\right)=
    \frac{1}{2\pi}\Log\,\Phi_A(2\pi,0)=\left[%
    \begin{array}{cc}
      -1 & 0 \\
      \frac{1}{2} & 0 \\
    \end{array}%
    \right],%
    \end{equation}
which is constant (as expected) since $\m\equiv 0$ on $\R$. Hence,
    $$
    e_R(t,0)=e^{Rt}=\left[%
    \begin{array}{cc}
      e^{-t} & 0 \\
     \frac{1}{2}-\frac{e^{-t}}{2} & 1 \\
    \end{array}%
    \right]%
    \qquad\text{and thus}\qquad e^{-Rt}=
    \left[%
    \begin{array}{cc}
      e^{t} & 0 \\
      \frac{1}{2}-\frac{e^{t}}{2} & 1 \\
    \end{array}%
    \right].
    $$

Furthermore, the $2\pi$-periodic matrix
$L(t)$ is given by
    \begin{align}
    L(t)&=\Phi_A(t,0)e_{R}^{-1}(t,0)\notag\\
    &=\Phi_A(t,0)e^{-Rt}\notag\\
    &=\left[%
    \begin{array}{cc}
      e^{-t} & 0 \\
      \frac{1}{2}-\frac{e^{-t}(\cos t+\sin t)}{2} & 1 \\
    \end{array}%
    \right]%
    \left[%
    \begin{array}{cc}
      e^{t} & 0 \\
      \frac{1}{2}-\frac{e^{t}}{2} & 1 \\
    \end{array}%
    \right]\notag\\
    &=\left[%
    \begin{array}{cc}
      1 & 0 \\
      \frac{1}{2} -\frac{\cos t+\sin t}{2} & 1 \\
    \end{array}%
    \right].\label{conts_example_L}%
    \end{align}
Using $R$ and $L(t)$ above, it is straightforward to verify that indeed \eqref{phi_exp_eq}, and equivalently, \eqref{Fdecomp} both hold here.

\subsection{Hybrid Example}\label{hybrid_example_1}

This example highlights the scope and generality of Theorem~\ref{Floquet2} on domains with nonconstant graininess.

Consider the time scale $\T=\mathbb{P}_{1,1}$ which
has a period of 2 and the regressive (on $\mathbb{P}_{1,1}$) time
varying matrix
    $$
    A(t)=\left[%
    \begin{array}{cc}
    g(t) & 1 \\
    0 & -3 \\
    \end{array}%
    \right],\qquad g(t)=-3+\sin(2\pi t),
    $$
with period 1. The transition matrix for the homogeneous periodic linear system of
dynamic equations
    \begin{equation}\label{hybrid_example_1_dynamics}
    X^\D(t)=
    \left[%
    \begin{array}{cc}
      g(t) & 1 \\
      0 & -3 \\
    \end{array}%
    \right]X(t),
    \end{equation}
is given by
    \begin{equation*}
    \Phi_A(t,0)=
    \left[%
    \begin{array}{cc}
      e_g(t,0) & C(t) \\
      0 & e_{-3}(t,0) \\
    \end{array}%
    \right],\quad\text{where}\quad C(t):=\int_0^t e_g(t,\s(\t))e_{-3}(\t,0)\D\t.
    \end{equation*}
Then $R(t)$ in the Floquet decomposition is given by
    \begin{equation}\label{hybrid_example_R}
    \begin{aligned}
    R(t)&:=\frac{1}{\m(t)}\left(\Phi_A(2,0)^\frac{\m(t)}{2}-I\right)\\
    &=\frac{1}{\m(t)}\left(\left[%
    \begin{array}{cc}
      -2e^{-3} & C(2) \\
      0 & -2e^{-3} \\
    \end{array}%
    \right]^{\frac{\m(t)}{2}}-\left[%
    \begin{array}{cc}
      1 & 0 \\
      0 & 1 \\
    \end{array}%
    \right]\right),
    \end{aligned}
    \end{equation}
which is nonconstant (as expected). Hence,
    $$
    e_R(t,0)=\left[%
    \begin{array}{cc}
      -2e^{-3} & C(2) \\
      0 & -2e^{-3} \\
    \end{array}%
    \right]^{\frac{t}{2}},
    $$
and thus,
    \begin{align*}
    e_R^{-1}(t,0)=\left[%
    \begin{array}{cc}
      -2e^{-3} & C(2) \\
      0 & -2e^{-3} \\
    \end{array}%
    \right]^{-\frac{t}{2}}.\\
    \end{align*}

Furthermore, the $2$-periodic matrix $L(t)$
is given by
    \begin{align}
    L(t)&=\Phi_A(t,0)e_{R}^{-1}(t,0)\notag\\
    &=\left[%
    \begin{array}{cc}
      e_g(t,0) & C(t) \\
      0 & e_{-3}(t,0) \\
    \end{array}%
    \right]\cdot\left[\begin{array}{cc}
      -2e^{-3} & C(2) \\
      0 & -2e^{-3} \\
    \end{array}%
    \right]^{-\frac{t}{2}}\label{hybrid_example_L}
    \end{align}
which is obviously $p$-periodic on $\P_{1,1}$.  Thus, the Floquet
decomposition of the transition matrix is
$\Phi_A(t,0)=L(t)e_{R}(t,0)$.

\section{Floquet Multipliers and Floquet Exponents}\label{FloquetMultipliersSection}

Suppose that $\Phi_A(t,t_0)$ is the transition matrix and
$\Phi(t)$ is the fundamental matrix at $t=\t$ (i.e. $\Phi(\t)=I$)
for the system \eqref{F}.  Then we can write any fundamental
matrix $\Psi(t)$ as
$$
\Psi(t)=\Phi(t)\Psi(\t)\quad\text{or}\quad\Psi(t)=\Phi_A(t,t_0)\Psi(t_0).
$$

\begin{definition}
Let $x_0\in\R^n$ be a nonzero vector and $\Psi(t)$ be any
fundamental matrix for the system \eqref{F}.  The vector solution
of the system with initial condition $x(t_0)=x_0$ is given by
$\Phi_A(t,t_0)x_0$. The operator $M:\R^n\to\R^n$ given by
    $$
    M(x_0):=\Phi_A(t_0+p,t_0)x_0=\Psi(t_0+p)\Psi^{-1}(t_0)x_0,
    $$
is called a {\it monodromy operator}. The eigenvalues of the
monodromy operator are called the {\it Floquet} (or {\it characteristic}) {\it multipliers} of the system \eqref{F}.
\end{definition}

The following theorem establishes that characteristic multipliers
are nonzero complex numbers intrinsic to the periodic system---they
do not depend on the choice of the fundamental matrix. This result
is analogous to the theorem dealing with the eigenvalues and
invertibility of monodromy operators in \cite{Chicone}, which can
also be referenced for proof.
\begin{theorem}
The following statements hold for the system \eqref{F}.
\begin{itemize}
    \item[(1)] Every monodromy operator is invertible.  In particular,
    every characteristic multiplier is nonzero.
    \item[(2)] If $M_1$ and $M_2$ are monodromy operators, then they
    have the same eigenvalues.  In particular, there are exactly
    $n$ characteristic multipliers, counting multiplicities.
\end{itemize}
\end{theorem}
%
%
%

With the Floquet normal form
$\Phi_A(t,t_0)=\Psi_1(t)\Psi_1^{-1}(t_0)=L(t)e_R(t,t_0)L^{-1}(t_0)$
of the transition matrix for the system \eqref{F} on one hand, and
the monodromy operator representation
    $$
    M(x_0)=\Phi_A(t_0+p,t_0)x_0=\Psi_1(t_0+p)\Psi_1^{-1}(t_0)x_0,
    $$
on the other, together we conclude
    $$
    \Phi_A(t_0+p,t_0)=\Psi_1(t_0+p)\Psi_1^{-1}(t_0)=L(t_0)e_R(t_0+p,t_0)L^{-1}(t
    _0).
    $$
Thus, the Floquet (or characteristic) multipliers of the system are the
eigenvalues of the matrix $e_R(t_0+p,t_0)$.  The (possibly complex)
scalar function $\g(t)$ is a {\it Floquet} ({\it or characteristic})
{\it exponent} of the $p$-periodic system \eqref{F} if $\l$ is a
Floquet multiplier and $e_\g(t_0+p,t_0)=\l$.

Theorem~\ref{spectralmap} will help us answer the question of
whether or not the eigenvalues of the matrix $R(t)$ in the
Floquet decomposition $\Phi_A(t,t_0)=L(t)e_R(t,t_0)$ are in fact
Floquet exponents.
\begin{theorem}[{\bf Spectral Mapping Theorem for Time
Scales}]\label{spectralmap}
Suppose that $R(t)$ is an $n\times n$
matrix as in Theorem~\ref{defofmatrixRthm}, with eigenvalues
$\g_1(t),...,\g_n(t)$, repeated according to multiplicities. Then
$\g_1^k(t),...,\g_n^k(t)$ are the eigenvalues of $R^k(t)$ and the
eigenvalues of $e_R$ are $e_{\g_1},...,e_{\g_n}$.
\end{theorem}

\begin{proof}
By induction for the dimension $n$, we start by stating that the
theorem is valid for $1\times 1$ matrices. Suppose that it is true
for all $(n-1)\times (n-1)$ matrices. For each fixed $t\in\T$, take
$\g_1(t)$ and let $v\neq 0$ denote a corresponding eigenvector such
that $R(t)v=\g_1(t)v$. Let ${\bf e}_1,...,{\bf e}_n$ denote the
usual basis of $\C^n$. There exists a nonsingular matrix $S$ such
that $Sv={\bf e}_1$. Thus we have $SR(t)S^{-1}{\bf e}_1=\g_1(t){\bf
e}_1,$ and the matrix $SR(t)S^{-1}$ has the block form
    $$
    SR(t)S^{-1}=\left[%
    \begin{array}{cc}
      \g_1(t) & * \\
      0 & \widetilde{R}(t) \\
    \end{array}%
    \right].
    $$
The matrix $SR^k(t)S^{-1}$ has the same block form, only with block
diagonal elements $\g_1^k(t)$ and $\widetilde{R}^k(t)$.  Clearly,
the eigenvalues of this block matrix are $\g_1^k(t)$ together with
the eigenvalues of $\widetilde{R}^k(t)$.  By induction, the
eigenvalues of $\widetilde{R}^k(t)$ are the $k$th powers of the
eigenvalues of $\widetilde{R}(t)$. This proves the first statement
of the theorem.

Since we chose the matrix $S$ so that $SR(t)S^{-1}$ is block
diagonal, by construction, the matrix $e_{SRS^{-1}}$ is also block
diagonal. We see that $e_{SRS^{-1}}$ has block diagonal form, with
block diagonal elements $e_{\g_1}$ and $e_{\tilde{R}}$. Using
induction, it follows that the eigenvalues of $e_{\tilde{R}}$ are
$e_{\g_2},...,e_{\g_n}$. Thus, the eigenvalues of
$e_{SRS^{-1}}=Se_RS^{-1}$ are $e_{\g_1},...,e_{\g_n}$.
\end{proof}
We now know that the eigenvalues of the matrix $e_R(t_0+p,t_0)$ are the Floquet multipliers and the eigenvalues of $R(t)$ are Floquet exponents.  However, in Theorem~\ref{nonuniquefloqexp}, we will see that although the Floquet exponents are the eigenvalues of the matrix $R(t)$, they are not unique.  We first introduce the definition of a purely imaginary number on an arbitrary time scale.

\begin{definition}Let $-\frac{\pi}{h}<\omega\leq\frac{\pi}{h}$.
The {\it Hilger purely imaginary number} $\i\omega$ is defined by
$
\i\omega=\frac{e^{i\omega h}-1}{h}.
$
For $z\in\C_h$, we have that $\i\text{Im}_h(z)\in\I_h$.  Also,
when $h=0$, $\i\omega=i\omega$.
\end{definition}

\begin{theorem}[{\bf Nonuniqueness of Floquet Exponents}]\label{nonuniquefloqexp}
Suppose $\g(t)\in\mathcal{R}$ is a \textup{(}possibly complex\textup{)} Floquet exponent, $\l$ is the corresponding characteristic multiplier of the $p$-periodic system \eqref{F} such that $e_{\g}(t_0+p,t_0)=\l$, and $\T$ is a $p$-periodic time scale. Then $\g(t)\oplus\i\frac{2\pi k }{p}$ is also a Floquet exponent
for all $k\in\Z$.
\end{theorem}

\begin{proof}
Observe that for any $k\in\Z$ and $t_0\in\T$,
    \begin{align*}
    e_{\g\oplus\i\frac{2\pi k
    }{p}}(t_0+p,t_0)&=e_{\g}(t_0+p,t_0)e_{\i\frac{2\pi
    k }{p}}(t_0+p,t_0)\\
    &=e_{\g}(t_0+p,t_0)\exp\left(\int_{t_0}^{t_0+p}\frac{\textup{Log}(1+\mu(\t)\i\frac{2\pi
    k}{p})}{\mu(\t)} \D\t\right)\\
    &=e_{\g}(t_0+p,t_0)\exp\left(\int_{t_0}^{t_0+p}\frac{\textup{Log}(1+\mu(\t)
    \frac{e^{i2\pi
    k\mu(\t)/p}-1}{\mu(\t)})}{\mu(\t)} \D\t\right)\\
    &=e_{\g}(t_0+p,t_0)\exp\left(\int_{t_0}^{t_0+p}\frac{\textup{Log}(e^{i2\pi
    k\mu(\t)/p})}{\mu(\t)} \D\t\right)\\
    &=e_{\g}(t_0+p,t_0)\exp\left(\int_{t_0}^{t_0+p}\frac{i2\pi
    k\mu(\t)/p}{\mu(\t)} \D\t\right)\\
    &=e_{\g}(t_0+p,t_0)\exp\left(\int_{t_0}^{t_0+p}\frac{i2\pi
    k}{p} \D\t\right)\\
    &=e_{\g}(t_0+p,t_0)e^{i2\pi k}\\
    &=e_{\g}(t_0+p,t_0).
    \end{align*}
\end{proof}



%

The next lemma will illustrate the periodic nature of a special exponential function which will be used in proving the nonuniqueness of Floquet exponents in Theorem~\ref{evalofR}
\begin{lemma}\label{circleiperiodic}
Let $\T$ be a $p$-periodic time scale and $k\in\Z$. Then the
functions $e_{\i\frac{2\pi k }{p}}(t,t_0)$ and
$e_{\ominus\i\frac{2\pi k }{p}}(t,t_0)$ are $p$-periodic.
\end{lemma}
\begin{proof}
Let $t\in\T$.  Then
    \begin{align*}
    e_{\i\frac{2\pi k }{p}}(t+p,t_0)
    &=\exp\({\frac{i2\pi k(t+p-t_0)}{p}}\)\\
    &=\exp\({\frac{i2\pi k(t-t_0)}{p}}\)\exp\({\frac{i2\pi kp}{p}}\)\\
    &=\exp\({\frac{i2\pi k(t-t_0)}{p}}\)\\
    &=e_{\i\frac{2\pi k }{p}}(t,t_0).
    \end{align*}
Therefore, $e_{\i\frac{2\pi k }{p}}(t,t_0)$ is a $p$-periodic
function.  The fact that $e_{\ominus\i\frac{2\pi k }{p}}(t,t_0)$
is $p$-periodic follows easily from the identity \cite{BoPe}
    $$
    e_{\ominus\i\frac{2\pi k }{p}}(t,t_0)=\frac{1}{e_{\i\frac{2\pi k
}{p}}(t,t_0)}.
    $$
\end{proof}

We now show that for any Floquet exponent, there exists a Floquet decomposition such that the Floquet exponent if an eigenvalue of the associated matrix $R(t)$.

\begin{theorem}[{\bf Special Floquet Decompositions}]\label{evalofR}
If $\g(t)$ is a Floquet
exponent for the system \eqref{F} and $\Phi_A(t,t_0)$ is the associated
transition matrix, then there exists a Floquet decomposition of the form
$\Phi_A(t,t_0)=L(t)e_R(t,t_0)$ where $\g(t)$ is an eigenvalue of
$R(t)$.
\end{theorem}

\begin{proof}
Consider the Floquet decomposition $\Phi_A(t,t_0)=\tilde{L}(t)e_{\tilde{R}}(t,t_0)$. By definition of the characteristic exponents, there is a characteristic multiplier $\l$ such that
$\l=e_\g(p+t_0,t_0)$, and, by Theorem~\ref{spectralmap}, there is
an eigenvalue $\tilde{\g}(t)$ of $\tilde{R}(t)$ such that
$e_{\tilde{\g}}(p+t_0,t_0)=\l$. Also, by Theorem~\ref{nonuniquefloqexp},
there is some integer $k$ such that $\tilde{\g}(t)=\g(t)\oplus
\i\frac{2\pi k}{p}$.

Set
    $$
    \begin{aligned}
    R(t)&:=\tilde{R}(t)\ominus \i\frac{2\pi k}{p}I,\\
    L(t)&:=\tilde{L}(t)e_{\i\frac{2\pi k}{p}I}(t,t_0).
    \end{aligned}
    $$
By this definition it is implied that $\tilde{R}(t)=R(t)\oplus \i\frac{2\pi k}{p}I$. Then $\g(t)$ is an eigenvalue of $R(t)$, $L(t)$ is a $p$-periodic function, and
    $$
    L(t)e_R(t,t_0)=\tilde{L}(t)e_{\i\frac{2\pi
    k}{p}I}(t,t_0)e_{R}(t,t_0)=\tilde{L}(t)e_{\i\frac{2\pi k}{p}I\oplus
    R}(t,t_0)=\tilde{L}(t)e_{\tilde{R}}(t,t_0).
    $$
It follows that $\Phi_A(t,t_0)=L(t)e_R(t,t_0)$ is another Floquet
decomposition where $\g(t)$ is an eigenvalue of $R(t)$.
\end{proof}

The following theorem classifies the types of
solutions that can arise with periodic systems.

\begin{theorem}\label{nontrivsol}
If $\l$ is a characteristic multiplier of the $p$-periodic system
\eqref{F} and $e_\g(t_0+p,t_0)=\l$ for some $t_0\in\T$, then there
exists a \textup{(}possibly complex\textup{)} nontrivial solution
of the form
    $$
    x(t)=e_\g(t,t_0)q(t)
    $$
where $q$ is a $p$-periodic function.  Moreover, for this solution
$x(t+p)=\l x(t)$.
\end{theorem}

\begin{proof}
Let $\Phi_A(t,t_0)$ be the transition matrix for \eqref{F}.  By
Theorem~\ref{evalofR}, there is a Floquet decomposition
$\Phi_A(t,t_0)=L(t)e_R(t,t_0)$ such that $\g(t)$ is an
eigenvalue of $R(t)$. There exists a vector $v\neq 0$ such that
$R(t)v=\g(t) v$. It follows that $e_R(t,t_0)v=e_\g(t,t_0)v$, and
therefore the solution $x(t):=\Phi_A(t,t_0)v$ can be represented
in the form
    $$
    x(t)=L(t)e_R(t,t_0)v=e_\g(t,t_0)L(t)v.
    $$
The solution required by the first part of the theorem is
obtained by defining $q(t):=L(t)v$.  The second part of
the theorem follows from
    \begin{align*}
    x(t+p)&=e_\g(t+p,t_0)q(t+p)\\
    &=e_\g(t+p,t_0+p)e_\g(t_0+p,t_0)q(t)\\
    &=e_\g(t_0+p,t_0)e_\g(t+p,t_0+p)q(t)\\
    &=e_\g(t_0+p,t_0)e_\g(t,t_0)L(t)v\\
    &=e_\g(t_0+p,t_0)x(t)\\
    &=\l x(t).
    \end{align*}
\end{proof}


The next result is motivated by \cite[Theorem 2.53]{Chicone} and is a direct consequence of Theorem~\ref{nontrivsol}.

\begin{corollary}\label{MultiplierCorollary}
Suppose that $\l_1,\dots,\l_n$ are the Floquet multipliers for the
$p$-periodic system \eqref{F}.
    \begin{itemize}
        \item[\textup{(1)}] If all the Floquet multipliers have modulus less than
        one, then the system \eqref{F} is exponentially stable.
        \item[\textup{(2)}] If all of the Floquet multipliers have modulus less
        than or equal to one, then the system \eqref{F} is stable.
        \item[\textup{(3)}] If at least one of the Floquet multipliers have
        modulus greater than one, then the system \eqref{F} is unstable.
    \end{itemize}
\end{corollary}

Theorem~\ref{nontrivsol} showed that if $\g(t)$ is a Floquet
exponent of \eqref{F}, then we can construct a nontrivial solution
of the form $x(t)=e_\g(t,t_0)q(t)$, where $q(t)$ is $p$-periodic. In
the next theorem it is shown that if two characteristic multipliers
$\l_1$ and $\l_2$ of the system \eqref{F} are distinct, then as in
Theorem~\ref{nontrivsol}, we can construct linearly independent
solutions $x_1$ and $x_2$ of \eqref{F}.

\begin{theorem}
Suppose that $\l_1,\l_2$ are characteristic multipliers of the
$p$-periodic system \eqref{F} and $\g_1(t),\g_2(t)$ are
Floquet exponents such that $e_{\g_i}(t_0+p,t_0)=\l_i$, $i=1,2$.  If $\l_1\neq\l_2$, then there exist
$p$-periodic functions $q_1(t),q_2(t)$ such that
    $$
    x_1(t)=e_{\g_1}(t,t_0)q_1(t)\quad\text{and}\quad
    x_2(t)=e_{\g_2}(t,t_0)q_2(t)
    $$
are linearly independent solutions of \eqref{F}.
\end{theorem}

\begin{proof}
As in Theorem~\ref{evalofR}, let $\Phi_A(t,t_0)=L(t)e_R(t,t_0)$ be
such that $\g_1(t)$ is an eigenvalue of $R(t)$ with corresponding
(nonzero) eigenvector $v_1$. Since $\l_2$ is an eigenvalue of the
monodromy matrix $\Phi_A(t_0+p,t_0)$, by Theorem~\ref{spectralmap}
there is an eigenvalue $\g(t)$ of $R(t)$ such that
$e_\g(t_0+p,t_0)=\l_2=e_{\g_2}(t_0+p,t_0)$. Hence
$\g_2(t)=\g(t)\oplus\i\frac{2\pi k}{p}$ for some $k\in\Z$.  Also,
$\g(t)\neq\g_1(t)$ since $\l_1\neq\l_2$.  Thus, if $v_2$ is a
nonzero eigenvector of $R(t)$ corresponding to the eigenvalue
$\g(t)$, then the eigenvectors $v_1$ and $v_2$ are linearly
independent.

As in the proof of Theorem~\ref{nontrivsol}, there are solutions of \eqref{F}
of the form
    $$
    x_1(t)=e_{\g_1}(t,t_0)L(t)v_1,\qquad
    x_2(t)=e_{\g}(t,t_0)L(t)v_2.
    $$
Because $x_1(t_0)=v_1$ and $x_2(t_0)=v_2$, these
solutions are linearly independent.  Finally, $x_2$ can be written
as
    $$
    x_2(t)=\left(e_{\g\oplus \i\frac{2\pi
    k}{p}}(t,t_0)\right)\left(e_{\ominus \i\frac{2\pi
    k}{p}}(t,t_0)L(t)v_2\right)=e_{\g_2}(t,t_0)\left(e_{\ominus
    \i\frac{2\pi k}{p}}(t,t_0)L(t)v_2\right),
    $$
where $q_2(t):=e_{\ominus \i\frac{2\pi k}{p}}(t,t_0)L(t)v_2$.
\end{proof}

\section{Examples Revisited}\label{ExamplesRevisited}
We now revisit the examples from Section~\ref{Examples}.  We show
for each of the three examples in Section~\ref{Examples} that the
original $R(t)$ matrices and the corresponding Floquet exponents
are not unique.  We conclude each example with a stability analysis
using the theorems developed from Section~\ref{FloquetTheory}.

\subsection{Discrete Time Example}

In Section~\ref{discrete_example_1}, the Floquet exponent we found for the system was
$\frac{\sqrt{3}}{2}-1$.  We show that
$\g:=-\frac{\sqrt{3}}{2}-1$ is also a Floquet exponent, but it
is not an eigenvalue of the original matrix $R$ in \eqref{discrete_example_R}.

Set $\tilde{R}:=R\ominus \i 2\pi Ik/p=R\ominus \i\pi I$, with $k=1$ and $p=2$,
as in Theorem~\ref{evalofR}. Then
    \begin{align*}
    \tilde{R}(t)&=\left[%
    \begin{array}{cc}
      \frac{\sqrt{3}}{2}-1 & 0 \\
      0 & \frac{\sqrt{3}}{2}-1 \\
    \end{array}%
    \right]\ominus
    \left[%
    \begin{array}{cc}
      \i\pi & 0 \\
      0 & \i\pi \\
    \end{array}%
    \right]\\
    &=\left[%
    \begin{array}{cc}
      \left(\frac{\sqrt{3}}{2}-1\right)\ominus\i\pi & 0 \\
      0 & \left(\frac{\sqrt{3}}{2}-1\right)\ominus\i\pi \\
    \end{array}%
    \right]\\
    &=\left[%
    \begin{array}{cc}
      -\frac{\sqrt{3}}{2}-1 & 0 \\
      0 & -\frac{\sqrt{3}}{2}-1 \\
    \end{array}%
    \right].
    \end{align*}
Hence,
    $$
    e_{\tilde{R}}(t,0)=(I+\tilde{R})^t=\left[%
    \begin{array}{cc}
      \left(-\frac{\sqrt{3}}{2}\right)^t & 0 \\
      0 & \left(-\frac{\sqrt{3}}{2}\right)^t \\
    \end{array}%
    \right],
    $$
and
    $$
    e_{\i\pi I}(t,0)=(-I)^t=\left[%
    \begin{array}{cc}
      (-1)^t & 0 \\
      0 & (-1)^t \\
    \end{array}%
    \right].
    $$

Next, using $L(t)$ from \eqref{discrete_example_L}, set
    $
    \tilde{L}(t):=L(t)e_{\i\pi I}(t,0)
    $
to obtain
    \begin{align*}
    \tilde{L}(t)&=\frac{1}{2}
    \left[%
    \begin{array}{cc}
      1+(-1)^t & \sqrt{3}+(-1)^t(-\sqrt{3})\\
      \sqrt{3}+(-1)^t(-\sqrt{3}) & 1+(-1)^t \\
    \end{array}%
    \right]\left[%
    \begin{array}{cc}
      (-1)^t & 0 \\
      0 & (-1)^t \\
    \end{array}%
    \right]\\
    &=\frac{(-1)^t}{2}
    \left[%
    \begin{array}{cc}
      1+(-1)^t & \sqrt{3}+(-1)^t(-\sqrt{3})\\
      \sqrt{3}+(-1)^t(-\sqrt{3}) & 1+(-1)^t \\
    \end{array}%
    \right].
    \end{align*}
Thus,
    $$
    \tilde{L}(t)e_{\tilde{R}}(t,0)=L(t)e_{\i\pi I}(t,0)e_{\tilde{R}}(t,0)
    =L(t)e_{\i\pi I\oplus \tilde{R}}(t,0)=L(t)e_{R}(t,0),
    $$
and so
    \begin{align*}
    L(t)e_{R}(t,0)&=\frac{1}{2}
    \left[%
    \begin{array}{cc}
      1+(-1)^t & \sqrt{3}+(-1)^t(-\sqrt{3})\\
      \sqrt{3}+(-1)^t(-\sqrt{3}) & 1+(-1)^t \\
    \end{array}%
    \right]\left[%
    \begin{array}{cc}
      \left(\frac{\sqrt{3}}{2}\right)^t & 0\\
      0 & \left(\frac{\sqrt{3}}{2}\right)^t \\
    \end{array}%
    \right]\\
    &=\frac{1}{2}
    \left[%
    \begin{array}{cc}
      1+(-1)^t & \sqrt{3}+(-1)^t(-\sqrt{3})\\
      \sqrt{3}+(-1)^t(-\sqrt{3}) & 1+(-1)^t \\
    \end{array}%
    \right]\left[%
    \begin{array}{cc}
      (-1)^t & 0 \\
      0 & (-1)^t \\
    \end{array}%
    \right]\left[%
    \begin{array}{cc}
      \left(-\frac{\sqrt{3}}{2}\right)^t & 0 \\
      0 & \left(-\frac{\sqrt{3}}{2}\right)^t \\
    \end{array}%
    \right]\\
    &=\frac{(-1)^t}{2}
    \left[%
    \begin{array}{cc}
      1+(-1)^t & \sqrt{3}+(-1)^t(-\sqrt{3})\\
      \sqrt{3}+(-1)^t(-\sqrt{3}) & 1+(-1)^t \\
    \end{array}%
    \right]\left[%
    \begin{array}{cc}
      \left(-\frac{\sqrt{3}}{2}\right)^t & 0 \\
      0 & \left(-\frac{\sqrt{3}}{2}\right)^t \\
    \end{array}%
    \right]\\
    &=\tilde{L}(t)e_{\tilde{R}}(t,0).
    \end{align*}
Therefore, $\Phi_A(t,0)=\tilde{L}(t)e_{\tilde{R}}(t,0)$ is another (distinct) Floquet
decomposition of the transition matrix for $A$. Moreover,
    $$
    \g:=(\sqrt{3}/2-1)\ominus\i\pi =-\sqrt{3}/2-1,
    $$
is a Floquet exponent as well as an eigenvalue of $R$ which
corresponds to the Floquet multiplier $\l=\frac{3}{4}$; that is,
    $$
    e_{(\sqrt{3}/2-1)\ominus\i\pi}(2,0)=e_{-\sqrt{3}/2-1}(2,0)=3/4.
    $$

In light of Corollary~\ref{MultiplierCorollary}, solutions of \eqref{discrete_example_1_dynamics} are exponentially stable since the Floquet multipliers $\l_{1,2}=3/4$ have modulus less
than 1.
\subsection{Continuous Time Example}

From Section~\ref{conts_example_1}, consider $R$ from \eqref{conts_example_R}. Set $\tilde{R}:=R\ominus \i 2\pi Ik/p=R-i\pi I$, with $k=1$ and $p=2\pi$,as in Theorem~\ref{evalofR}. Then,
    $$
    \tilde{R}=\left[%
    \begin{array}{cc}
      -1 & 0 \\
      \frac{1}{2} & 0 \\
    \end{array}%
    \right]-
    \left[%
    \begin{array}{cc}
      i & 0 \\
      0 & i \\
    \end{array}%
    \right]
    =\left[%
    \begin{array}{cc}
      -1-i & 0 \\
      \frac{1}{2} & -i\\
    \end{array}%
    \right],
    $$
so that
    $$
    e_{\tilde{R}}(t,0)=e^{\tilde{R}t}=\left[%
    \begin{array}{cc}
      e^{(-1-i)t} & 0 \\
      \frac{e^{-it}}{2}-\frac{e^{(-1-i)t}}{2} & e^{-it} \\
    \end{array}%
    \right],
    $$
and
    $$\Phi_A(2\pi,0)=e^{2\pi \tilde{R}}=\left[%
    \begin{array}{cc}
      e^{-2\pi} & 0 \\
      \frac{1}{2}-\frac{e^{-2\pi}}{2} & 1 \\
    \end{array}%
    \right].
    $$

Next, using $L(t)$ from \eqref{conts_example_L}, set $\tilde{L}(t):=L(t)e_{\i I}(t,0)=L(t)e^{iIt}$
to obtain
    $$
    \tilde{L}(t)=\left[%
    \begin{array}{cc}
      1 & 0 \\
      \frac{1}{2} -\frac{\cos t+\sin t}{2} & 1 \\
    \end{array}%
    \right]\cdot\left[%
    \begin{array}{cc}
      e^{it} & 0 \\
      0& e^{it} \\
    \end{array}%
    \right]=\left[%
    \begin{array}{cc}
      e^{it} & 0 \\
      \frac{e^{it}}{2}-\frac{e^{it}(\cos t+\sin t)}{2}& e^{it} \\
    \end{array}%
    \right].
    $$
Thus,
    $$
    \tilde{L}(t)e_{\tilde{R}}(t,0)=L(t)e^{iIt}e^{\tilde{R}t}=L(t)e^{(iI+ \tilde{R})t}=L(t)e^{Rt},
    $$
and so
    \begin{align*}
    L(t)e^{Rt}&=\left[%
    \begin{array}{cc}
      1 & 0 \\
      \frac{1}{2} -\frac{\cos t+\sin t}{2} & 1 \\
    \end{array}%
    \right]\left[%
    \begin{array}{cc}
      e^{-t} & 0 \\
      \frac{1}{2}-\frac{e^{-t}}{2} & 1 \\
    \end{array}%
    \right]\\
    &=\left[%
    \begin{array}{cc}
      1 & 0 \\
      \frac{1}{2} -\frac{\cos t+\sin t}{2} & 1 \\
    \end{array}%
    \right]\left[%
    \begin{array}{cc}
      e^{it} & 0 \\
      0& e^{it} \\
    \end{array}%
    \right]\left[%
    \begin{array}{cc}
      e^{(-1-i)t} & 0 \\
      \frac{e^{-it}}{2}-\frac{e^{(-1-i)t}}{2}  & e^{-it} \\
    \end{array}%
    \right]\\
    &=\left[%
    \begin{array}{cc}
      e^{it} & 0 \\
      \frac{e^{it}}{2}-\frac{e^{it}(\cos t+\sin t)}{2}& e^{it} \\
    \end{array}%
    \right]\left[%
    \begin{array}{cc}
      e^{(-1-i)t} & 0 \\
      \frac{e^{-it}}{2}-\frac{e^{(-1-i)t}}{2}  & e^{-it} \\
    \end{array}%
    \right]\\
    &=\tilde{L}(t)e^{\tilde{R}t}.
    \end{align*}
Therefore, $\Phi_A(t,0)=\tilde{L}(t)e_{\tilde{R}}(t,0)$ is another (distinct) Floquet
decomposition of the transition matrix for $A$. Moreover, $\g_1:=-1-i$, $\g_2:=-i$
are Floquet exponents as well as eigenvalues of $R$ which correspond
to the Floquet multipliers $\l_1=e^{-2\pi}$, $\l_2=1$,
respectively. That is, $e^{-2\pi-2\pi i}=e^{-2\pi}$ and $e^{-2\pi
i}=1$.

In light of Corollary~\ref{MultiplierCorollary}, solutions of \eqref{conts_example_1_dynamics} are uniformly stable since the Floquet multipliers satisfy $|\l_i|\le 1$, $i=1,2$.

\subsection{Hybrid Example}

In Section~\ref{hybrid_example_1}, consider $R$ from \eqref{hybrid_example_R} which was given by
    \begin{align*}
    R(t)&=\frac{1}{\m(t)}\left(\left[%
    \begin{array}{cc}
      -2e^{-3} & C(2) \\
      0 & -2e^{-3} \\
    \end{array}%
    \right]^{\frac{\m(t)}{2}}-\left[%
    \begin{array}{cc}
      1 & 0 \\
      0 & 1 \\
    \end{array}%
    \right]\right),
    \end{align*}
where $-2e^{-3}=e_{-3}(2,0)$ on $\T$.  The eigenvalues of $R(t)$ are
    $$
    \g_1(t):=\lim_{s\searrow\mu(t)}s^{-1}((e_{-3}(2,0))^{-\frac{s}{2}}-1)
    \quad\text{ and }\quad
    \g_2(t):=\lim_{s\searrow\mu(t)}s^{-1}((e_{-3+\sin(2\pi t)}
    (2,0))^{-\frac{s}{2}}-1).
    $$
Set $\tilde{R}:=R\ominus \i 2\pi Ik/p=R\ominus \i\pi I$, with $k=1$ and $p=2$, as in Theorem~\ref{evalofR}. Then,
    \begin{align*}
    \tilde{R}(t)&=\frac{1}{\m(t)}\left(\left[%
    \begin{array}{cc}
      -2e^{-3} & C(2) \\
      0 & -2e^{-3} \\
    \end{array}%
    \right]^{\frac{\m(t)}{2}}-\left[%
    \begin{array}{cc}
      1 & 0 \\
      0 & 1 \\
    \end{array}%
    \right]\right)\ominus\left[%
    \begin{array}{cc}
      \i\pi & 0 \\
      0 & \i\pi \\
    \end{array}%
    \right],
    \end{align*}
and thus,
    \begin{align*}
    e_{\tilde{R}}(t,0)&=\left[%
    \begin{array}{cc}
      -2e^{-3} & C(2) \\
      0 & -2e^{-3} \\
    \end{array}%
    \right]^{\frac{t}{2}}e_{\ominus\i\pi I}(t,0)\\
    &=\left[%
    \begin{array}{cc}
      -2e^{-3} & C(2) \\
      0 & -2e^{-3} \\
    \end{array}%
    \right]^{\frac{t}{2}}e_{\i\pi I}(0,t)\\
    &=\left[%
    \begin{array}{cc}
      -2e^{-3} & C(2) \\
      0 & -2e^{-3} \\
    \end{array}%
    \right]^{\frac{t}{2}}\left[%
    \begin{array}{cc}
      e^{-i\pi t} & 0 \\
      0 & e^{-i\pi t} \\
    \end{array}%
    \right].\\
    \end{align*}

Next, using $L(t)$ from \eqref{hybrid_example_L}, set $\tilde{L}(t):=L(t)e_{\i\pi I}(t,0)$
to obtain
    \begin{align*}
    \tilde{L}(t)&=\left[%
    \begin{array}{cc}
      e_g(t,0) & C(t) \\
      0 & e_{-3}(t,0) \\
    \end{array}%
    \right]\left[\begin{array}{cc}
      -2e^{-3} & C(2) \\
      0 & -2e^{-3} \\
    \end{array}%
    \right]^{-\frac{t}{2}}\left[%
    \begin{array}{cc}
      e^{i\pi t} & 0 \\
      0 & e^{i\pi t} \\
    \end{array}%
    \right].
    \end{align*}
Thus,
    $$
    \tilde{L}(t)e_{\tilde{R}}(t,0)=L(t)e_{\i\pi I}(t,0)e_{\tilde{R}}(t,0)=L(t)e_{\tilde{R}}(t,0)e_{\i\pi
    I}(t,0)=L(t)e_{\tilde{R}\oplus\i\pi I}(t,0)=L(t)e_{R}(t,0),
    $$
and so
    \begin{align*}
    &L(t)e_{R}(t,0)=\left(\left[%
    \begin{array}{cc}
      e_g(t,0) & \int_0^t e_g(t,\s(\t))e_{-3}(\t,0)\D\t \\
      0 & e_{-3}(t,0) \\
    \end{array}%
    \right]\left[\begin{array}{cc}
      -2e^{-3} & \int_0^2 e_g(2,\s(\t))e_{-3}(\t,0)\D\t \\
      0 & -2e^{-3} \\
    \end{array}%
    \right]^{-\frac{t}{2}}\right)\\
    &\cdot\left[\begin{array}{cc}
      -2e^{-3} & \int_0^2 e_g(2,\s(\t))e_{-3}(\t,0)\D\t \\
      0 & -2e^{-3} \\
    \end{array}%
    \right]^{\frac{t}{2}}\\
    &=\left(\left[%
    \begin{array}{cc}
      e_g(t,0) & \int_0^t e_g(t,\s(\t))e_{-3}(\t,0)\D\t \\
      0 & e_{-3}(t,0) \\
    \end{array}%
    \right]\left[\begin{array}{cc}
      -2e^{-3} & \int_0^2 e_g(2,\s(\t))e_{-3}(\t,0)\D\t \\
      0 & -2e^{-3} \\
    \end{array}%
    \right]^{-\frac{t}{2}}\right)\\
    &\cdot\left(\left[%
    \begin{array}{cc}
      -2e^{-3} & \int_0^2 e_g(2,\s(\t))e_{-3}(\t,0)\D\t \\
      0 & -2e^{-3} \\
    \end{array}%
    \right]^{\frac{t}{2}}\left[%
    \begin{array}{cc}
      e^{-i\pi t} & 0 \\
      0 & e^{-i\pi t} \\
    \end{array}%
    \right]\right)\cdot\left[%
    \begin{array}{cc}
      e^{i\pi t} & 0 \\
      0 & e^{i\pi t} \\
    \end{array}%
    \right]\\
    &=\left(\left[%
    \begin{array}{cc}
      e_g(t,0) & \int_0^t e_g(t,\s(\t))e_{-3}(\t,0)\D\t \\
      0 & e_{-3}(t,0) \\
    \end{array}%
    \right]\left[\begin{array}{cc}
      -2e^{-3} & \int_0^2 e_g(2,\s(\t))e_{-3}(\t,0)\D\t \\
      0 & -2e^{-3} \\
    \end{array}%
    \right]^{-\frac{t}{2}}\cdot\left[%
    \begin{array}{cc}
      e^{i\pi t} & 0 \\
      0 & e^{i\pi t} \\
    \end{array}%
    \right]\right)\\
    &\cdot\left(\left[%
    \begin{array}{cc}
      -2e^{-3} & \int_0^2 e_g(2,\s(\t))e_{-3}(\t,0)\D\t \\
      0 & -2e^{-3} \\
    \end{array}%
    \right]^{\frac{t}{2}}\left[%
    \begin{array}{cc}
      e^{-i\pi t} & 0 \\
      0 & e^{-i\pi t} \\
    \end{array}%
    \right]\right)\\
    &=\tilde{L}(t)e_{\tilde{R}}(t,0).
    \end{align*}
Therefore, $\Phi_A(t,0)=\tilde{L}(t)e_{\tilde{R}}(t,0)$ is another (distinct) Floquet
decomposition of the transition matrix for $A$. Moreover, $\g(t):=-3\ominus\i\pi$ is
another Floquet exponent as well as an eigenvalue of $R(t)$ which
corresponds to the Floquet multipliers $\l_{1,2}=-2e^{-3}$. That is,
$e_{-3\ominus\i\pi}(2,0)=e_{-3}(2,0)=-2e^{-3}$.

In light of Corollary~\ref{MultiplierCorollary}, solutions of \eqref{hybrid_example_1_dynamics} are exponentially stable since the Floquet multipliers satisfy $|\l_i|< 1$, $i=1,2$.

%

\section{Floquet Theory, Stability, and Dynamic Eigenpairs}

Consider the $p$-periodic regressive system
    \begin{equation}\label{Aprob}
    x^\D(t)=A(t)x(t), \quad x(t_0)=x_0.
    \end{equation}
Recall that the $R$ matrix in the Floquet decomposition of $\Phi_A$ is given by
    \begin{equation}\label{Rdef}
    R(t):=\lim_{s\searrow\m(t)}\frac{1}{s}\left(\Phi_A(t_0+p,t_0)^\frac{s}{p}-I\right),
    \end{equation}
and consider the uniformly regressive system associated with \eqref{Aprob},
    \begin{equation}\label{Rprob}
    z^\D(t)=R(t)z(t), \quad z(t_0)=x_0.
    \end{equation}
From Theorem~\ref{relationship_thm}, solutions to these two problems are related via $z(t)=L^{-1}(t)x(t)$ where $L(t)$ is the Lyapunov transformation from \eqref{L}.

Given a constant $n\times n$ matrix $M$, let $C$ be the nonsingular matrix that transforms $M$ into its Jordan canonical form,
    $$
    J:=C^{-1}MC=\diag\left[J_{m_1}(\l_1),\dots,J_{m_k}(\l_k)\right],
    $$
where $k\leq n$, $\sum_{i=1}^k m_i=n$, $\l_i$ are the eigenvalues of $M$ (where some of them may be equal), and $J_m(\l)$ is a $m\times m$ Jordan block,
\begin{equation}\label{jordan block form}
J_m(\l)=
\left[
  \begin{array}{ccccc}
    \l & 1 &  &  &  \\
     & \l & 1 &  &  \\
     &   & \ddots & \ddots &  \\
     &   &   & \ddots & 1 \\
     &   &   &   & \l \\
  \end{array}
\right].
\end{equation}
We introduce a definition from~\cite{Poochy} which gives a bound on the eigenvalues of a system.  Then we state a lemma which proves that the system~\eqref{Rprob} associated with~\eqref{Aprob} via a Floquet decomposition is necessarily uniformly regressive.

\begin{definition}\cite{Poochy}
The scalar function $\g:\T\to\C$ is {\it uniformly regressive} if there exists a constant $\delta> 0$ such that
$0<\delta^{-1}\leq|1+\mu(t)\g(t)|$, for all $t\in\T^\kappa$.
\end{definition}


\begin{lemma}
Each eigenvalue of the matrix $R(t)$ in \eqref{Rprob} is uniformly regressive.
\end{lemma}

\begin{proof}
We show by direct substitution that the eigenvalues of $R(t)$
are uniformly regressive. Using Corollary~\ref{cor_evals of M
are the evals of R}, let
$\g_i(t):=\lim_{s\searrow\m(t)}\frac{\l_i^{\frac{s}{p}}-1}{s}$
be any of the $k\leq n$ distinct eigenvalues of $R(t)$. Recall
that in this $p$-periodic setting $0\leq\m(t)\leq p$.  Suppose
that $|\l_i|\geq 1$. Then,
    \begin{align*}
    |1+\mu(t)\g_i(t)|=\lim_{s\searrow\m(t)}\left|1+s\frac{\l_i^{\frac{s}{p}}-1}{s}\right|
    =\lim_{s\searrow\m(t)}|\l_i^{\frac{s}{p}}|\geq 1.
    \end{align*}
On the other hand, suppose $0<|\l_i|<1$.  Then,
    \begin{align*}
    |1+\mu(t)\g_i(t)|=\lim_{s\searrow\m(t)}\left|1+s\frac{\l_i^{\frac{s}{p}}-1}{s}\right|
    =\lim_{s\searrow\m(t)}|\l_i^{\frac{s}{p}}|\geq |\l_i|.
    \end{align*}
Thus $\delta^{-1}:=\min\{1,|\l_1|,\dots,|\l_k|\}$ will
suffice as the bound for uniform regressivity.
\end{proof}





The following definition will be used to obtain the main result of this section.

\begin{definition}\label{def Hilger Circle}
Let $\C_\m:=\{z\in\C:z\neq-\frac{1}{\m(t)}\}$. Given an element $t\in\T^\kappa$ with $\m(t)>0$, we define the {\it Hilger circle} as
    $$
    \H_t:=\{z\in\C_\m:\left|1+z\m(t)\right|<1\}.
    $$
Thus, we define the {\it closed Hilger circle} as
    $$
    \overline{\H}_t:=\{z\in\C_\m:\left|1+z\m(t)\right|\leq 1\},
    $$
Likewise, when $\m(t)=0$, we define the Hilger circle as
    $$
    \H_t:=\{z\in\C:\textup{Re}\,z<0\}
    $$
and the closed Hilger circle as
    $$
    \overline{\H}_t:=\{z\in\C:\textup{Re}\,z\leq 0\}.
    $$
\end{definition}



Motivated by the work of Wu \cite{Wu} on $\R$, in \cite{DaDaModal} we introduced the following.

\begin{definition}\label{Dyn Eval prob}
A nonzero, delta differentiable vector $w(t)$ is said to be a {\it
dynamic eigenvector} of $M(t)$ associated with the {\it dynamic eigenvalue} $\xi(t)$ if the pair satisfies the {\it dynamic eigenvalue problem}
    \begin{equation}\label{dyn_eval}
    w^\D(t)=M(t)w(t)-\xi(t)w^\s(t), \quad t\in\T^\kappa.
    \end{equation}
We call $\{\xi(t),w(t)\}$ a {\it dynamic eigenpair}.

The nonzero, delta differentiable vector
    \begin{equation}\label{modeVector}
    m_i(t):=e_{\xi_i}(t,t_0)w_i(t),
    \end{equation}
is the {\it mode vector} of $M(t)$ associated with
the dynamic eigenpair $\{\xi_i(t),w_i(t)\}$.
\end{definition}

The following lemma proves that given any regressive linear dynamic system (not necessarily periodic), there always exists a set of $n$ dynamic eigenpairs, with linearly independent dynamic eigenvectors.

\begin{lemma}\label{Lem Existence of Dyn Epairs}
Given the $n\times n$ regressive matrix $M$, there
always exists a set of $n$ dynamic eigenpairs with linearly independent dynamic eigenvectors.  Each of the eigenpairs satisfies the vector dynamic eigenvalue problem~\eqref{dyn_eval} associated with $M$.  Furthermore, when the $n$ vectors form the columns of $W(t)$, then $W(t)$ satisfies the equivalent matrix dynamic eigenvalue problem
    \begin{equation}\label{dyn_eval prob MATRIX FORM}
    W^\D(t)=M(t)W(t)-W^\s(t)\Xi(t),\quad\text{where}\quad\Xi(t):=\diag[\xi_1(t),\ldots,\xi_n(t)].
    \end{equation}
\end{lemma}

\begin{proof}
Let $\{\xi_i(t)\}_{i=1}^n$ be a set of (not necessarily
distinct) regressive functions.  Then the $n\times n$ nonsingular
matrix $W(t)$ defined by
    $$
    W(t):=\Phi_M(t,t_0)e_{\ominus\Xi}(t,t_0),
    $$
has as its columns the associated $n$ linearly independent
dynamic eigenvectors $\{w_i(t)\}_{i=1}^n$.  By direct substitution into~\eqref{dyn_eval prob MATRIX FORM}, the proof is complete.
\end{proof}

The next theorem is from \cite{DaDaModal}. We show that the stability of a linear dynamic system can be completely determined by the mode vectors $m_i$, $1\leq i\leq n$, as constructed in~\eqref{modeVector}.

\begin{theorem}[{\bf Stability Via Modal Vectors}]\label{mode vector stability theorem}
Solutions to the uniformly regressive \textup{(}but not necessarily periodic\textup{)} time varying
linear dynamic system~\eqref{A} are:
\begin{itemize}
  \item[(i)] stable if and only if there exists a $\g>0$ such that every mode vector $m_i(t)$ of $A(t)$ satisfies $||m_i(t)||\leq\g<\infty,\; t>t_0,$ for all $1\le i\le n$,
  \item[(ii)] asymptotically stable if and only if, in addition to \textup{(i)}, $||m_i(t)||\to 0,\; t>t_0,$ for
  all $1\le i\le n$,
  \item[(iii)] exponentially stable if and only if there exists $\g,\l>0$ with $-\l\in\mathcal{R^+}(\T,\R)$ such that $||m_i(t)||\leq\g e_\l(t,t_0),\; t>t_0,$ for all $1\le i\le n$.
\end{itemize}
\end{theorem}

\begin{proof}
Let $\{\xi_i(t),w_i(t)\}_{i=1}^n$ be a set of $n$ dynamic eigenpairs with linearly
independent dynamic eigenvectors associated with the system matrix
in~\eqref{A}. The transition matrix can be represented by
\begin{equation}\label{transition Matrix
def with Wn}
\Phi_A(t,t_0)=W(t)e_{\Xi}(t,t_0)W^{-1}(t_0)
\end{equation}
where $W(t):=\left[w_1(t),\;w_2(t),\cdots,w_n(t)\right]$ and
$\Xi(t):=\diag[\xi_1(t),\ldots,\xi_n(t)]$.

Denoting the matrix $W^{-1}(t_0)$ as
\begin{equation}\notag
W^{-1}(t_0):=\left[
               \begin{array}{c}
                 v_1^T(t_0) \\
                v_2^T(t_0) \\
                 \vdots \\
                 v_n^T(t_0) \\
               \end{array}
             \right],
\end{equation}
we now have the reciprocal basis of each $w_i(t_0)$ given by the
row vectors $v_i^T(t_0)$.

Because $\Xi(t)$ is a diagonal matrix,~\eqref{transition Matrix
def with Wn} can be rewritten as
    \begin{equation}\label{transition matrix Def with Fi}
    \Phi_A(t,t_0)=\sum_{i=1}^n e_{\xi_i}(t,t_0)W(t)F_iW^{-1}(t_0),
    \end{equation}
where $F_i:=\delta_{ij}$ is $n\times n$. Observing that
$v_i^T(t)w_j(t)=\delta_{ij}$ for all $t\in\T$, we may rewrite $F_i$ as
    \begin{equation}\label{Fi rewritten}
    F_i=W^{-1}(t)\left[0,\dots,0,\;w_i(t),0,\ldots,0\right].
    \end{equation}
Substituting~\eqref{Fi rewritten} into~\eqref{transition matrix
Def with Fi} yields
    \begin{equation}\label{transition matrix Def with Fi Subs}
    \Phi_A(t,t_0)=\sum_{i=1}^n
    e_{\xi_i}(t,t_0)w_i(t)v_i^T(t_0)=\sum_{i=1}^n m_i(t)v_i(t_0).
    \end{equation}
The proof concludes easily from~\eqref{transition matrix Def
with Fi Subs}.
\end{proof}

We now set the stage for the main result of the section.  In the next theorem, we will show that given the system matrix $R(t)$ from~\eqref{Rprob}, we can choose a set of $n$ linearly independent dynamic eigenvectors that have a growth rate that is bounded by a finite sum of generalized polynomials.
\begin{theorem}\label{thm-bd on dyn evecs}
Given the set of traditional eigenvalues $\{\g_i(t)\}_{i=1}^n$ from the matrix $R(t)$ in \eqref{Rdef}, let $\{w_i(t)\}_{i=1}^n$ denote the corresponding linearly independent dynamic eigenvectors as defined by Lemma~\ref{Lem Existence of Dyn Epairs}. Then $\{\g_i(t),w_i(t)\}_{i=1}^n$ is a set of dynamic eigenpairs of $R(t)$ with the property that each $w_i(t)$ is bounded by at most a finite sum of constant multiples of generalized polynomials.  In other words, there exists positive constants $D_i>0$ such that
\begin{equation}\label{dyn evec bound}
||w_i(t)||\leq D_i\sum_{k=0}^{m_i-1}h_k(t,t_0),
\end{equation}
where $m_i$ is the dimension of the Jordan block which contains the $i$th eigenvalue, for all $1\leq i\leq n$.
\end{theorem}

\begin{proof}
The fact that  $\{\g_i(t),w_i(t)\}_{i=1}^n$ is a set of dynamic eigenpairs of $R(t)$ follows immediately from Lemma~\ref{Lem Existence of Dyn Epairs}. We now show each dynamic eigenvector $w_i$ satisfies the bound in~\eqref{dyn evec bound}.  For an appropriately chosen nonsingular $n\times n$ constant matrix $C$, we define the Jordan form of the constant matrix $\Phi_A(t_0+p,t_0)$ by
\begin{align}\label{Jordan Decomp of Phi(t_0+p,t_0)}
J&:=C^{-1}\Phi_A(t_0+p,t_0)C\\
&=\notag \left[
    \begin{array}{cccc}
      J_{m_1}(\l_1) &  &  &  \\
       & J_{m_2}(\l_2) &  &  \\
       &  & \ddots &  \\
       &  &  & J_{m_d}(\l_d) \\
    \end{array}
  \right]_{n\times n,}
\end{align}
where $d\leq n$, $\sum_{i=1}^d m_i=n$, $\l_i$ are the eigenvalues of $\Phi_A(t_0+p,t_0)$ (where some of them may be equal), and $J_m(\l)$ is a $m\times m$ Jordan block as defined in~\eqref{jordan block form}.

Using the matrix $C$ from~\eqref{Jordan Decomp of Phi(t_0+p,t_0)}, define
    \begin{align}
    K(t):=&\,C^{-1}R(t)C\label{doh}\\
    =&\,C^{-1}\lim_{s\searrow\m(t)}\frac{1}{s}\left(\Phi_A(t_0+p,t_0)^\frac{s}{p}-I\right)C\notag \\
    =&\,\lim_{s\searrow\m(t)}\frac{1}{s}\left(C^{-1}\Phi_A(t_0+p,t_0)^\frac{s}{p}C-I\right)\notag\\
    =&\,\lim_{s\searrow\m(t)}\frac{1}{s}\left(J^\frac{s}{p}-I\right),\notag
    \end{align}
where the penultimate equality is via Theorem~\ref{Thm-real power of M relationship with real power of JCF}.  Note that $K(t)$ has the block diagonal form
\begin{equation}\label{form of K(t)}
K(t)=\left[
    \begin{array}{cccc}
      K_{1}(t) &  &  &  \\
       & K_{2}(t) &  &  \\
       &  & \ddots &  \\
       &  &  & K_{d}(t) \\
    \end{array}
  \right]_{n\times n,},
\end{equation}
where each $K_i(t)$ from~\eqref{form of K(t)} has the form
$$
K_i(t):=
    \left[\begin{array}{ccccc}
    \frac{\l_i ^{\frac{\mu(t)}{p}}-1}{\mu(t)} & \frac{\l_i ^{\frac{\mu(t)}{p}-1}}{1! p} &\frac{(\frac{\mu(t)}{p}-1)\l_i ^{\frac{\mu(t)}{p}-2}}{2! p}  & \ldots & \frac{(\frac{\mu(t)}{p}-1)\cdots(\frac{\mu(t)}{p}-(n-2))\l_i ^{\frac{\mu(t)}{p}-(n-1)}}{(n-1)! p} \\
    & \frac{\l_i ^{\frac{\mu(t)}{p}}-1}{\mu(t)} & \frac{\l_i ^{\frac{\mu(t)}{p}-1}}{1! p} & \ldots & \frac{(\frac{\mu(t)}{p}-1)\cdots(\frac{\mu(t)}{p}-(n-3))\l_i ^{\frac{\mu(t)}{p}-(n-2)}}{(n-2)! p} \\
    &  & \frac{\l_i ^{\frac{\mu(t)}{p}}-1}{\mu(t)}  & \ddots & \vdots \\
    & &    &  \ddots & \frac{\l_i ^{\frac{\mu(t)}{p}-1}}{1! p} \\
    & &    & & \frac{\l_i ^{\frac{\mu(t)}{p}}-1}{\mu(t)} \\
    \end{array}%
    \right]_{{m_i\times m_i}.}
$$

By hypothesis, the dynamic eigenvalues of $R(t)$ have been defined as its traditional eigenvalues $\g_i(t)= \lim_{s\searrow\mu(t)}\frac{\l_i ^{\frac{s}{p}}-1}{s}$, for $1\leq i\leq n$.  Note that since $R(t)$ and $K(t)$ are similar, they have the same traditional eigenvalues, with corresponding multiplicities.  Furthermore, choosing the $n$ dynamic eigenvalues of $K(t)$ to be its $n$ traditional eigenvalues, we claim the corresponding dynamic eigenvectors $\{u_i(t)\}_{i=1}^n$ are defined to be $u_i(t):=C^{-1}w_i(t)$.

To prove this, since $\{\g_i(t),w_i(t)\}_{i=1}^n$ is a set of dynamic eigenpairs for $R(t)$, by definition we have
$$
w_i^\D(t)=R(t)w_i(t)-\g_i(t)w_i^\s(t),
$$
for all $1\leq i\leq n$.  We can now show that $\{\g_i(t),u_i(t)\}_{i=1}^n$ is a set of dynamic eigenpairs for $K(t)$:
\begin{align}\label{dyn eval prob for nu_i}
u_i^\D(t)&=C^{-1}w_i^\D(t)\\\notag
&=C^{-1}R(t)w_i(t)-C^{-1}\g_i(t)w_i^\s(t)\\\notag
&=K(t)C^{-1}w_i(t)-\g_i(t)C^{-1}w_i^\s(t)\\\notag
&=K(t)u_i(t)-\g_i(t)u_i^\s(t),\notag
\end{align}
for all $1\leq i\leq n$.

We now show that each of the $u_i(t)$ is bounded by a finite sum of constant multiples of generalized polynomials.

Since each of the $n$ dynamic eigenpairs $\{\g_i(t),u_i(t)\}_{i=1}^n$ satisfy the dynamic eigenvalue problem~\eqref{dyn eval prob for nu_i}, we can now find the structure of each of the dynamic eigenvectors.  Choose the $i$th block of $K(t)$, with dimension $m_i\times m_i$.  We must then solve the $m_i\times m_i$ linear dynamic system
    \begin{align}\label{form of K}
    \nu^\D(t)=\tilde{K}_i(t)\nu(t)=\left[%
    \begin{array}{ccccc}
      0& \frac{1}{\l_i p} &\frac{(\m(t)-p)}{2\l_i^2 p^2}  & \ldots & \frac{(\m(t)-p)\cdots(\m(t)-p(n-2))}{(n-1)!\;\l_i^{n-1} p^{n-1}} \\
       & 0 & \frac{1}{\l_i p} & \ldots & \frac{(\m(t)-p)\cdots(\m(t)-p(n-3))}{(n-2)!\;\l_i^{n-2} p^{n-2}}\\
       &  &   & \ddots & \vdots \\
       & &    & \ddots & \vdots \\
       & &    &  & \frac{1}{\l_i p} \\
       & &    & & 0
    \end{array}%
    \right]\nu(t),
    \end{align}
where $\tilde{K}_i(t):=K_i(t)\ominus\g_i(t)I$.  Since there are $m_i$ linearly independent solutions
to~\eqref{form of K}, we will denote each of the $m_i\times 1$ solutions by $\nu_{i,j}(t)$, where $i$ corresponds to the $i$th block matrix $K_i(t)$ of $K(t)$ and $j=1,\ldots,m_i$.  For $1\leq i\leq d$, define $l_i:= \sum_{s=0}^{i-1}m_s$, with $m_0:=0$.  The form for an arbitrary $n\times 1$ column vector $u_{l_i+j}$, with $1\leq j\leq m_i$, is
\begin{equation}\label{vec form u_i}
u_{l_i+j}^T(t)=[\underbrace{0,\ldots,0,\ldots,0}_{m_1+\dots+m_{i-1}},\underbrace{\nu_{i,j}^T(t)}_{m_i},
\underbrace{0,\ldots,0,\ldots,0}_{m_{i+1}+\dots+m_d} ]_{1\times n}.
\end{equation}
By combining all $n$ vector solutions from~\eqref{dyn eval prob for nu_i}, the solution to the equivalent $n\times n$ matrix dynamic equation
$$
U^\D(t)=K(t)U(t)-U^\s(t)\Gamma(t),
$$
where $\Gamma(t):=\diag\left[\g_1(t),\ldots,\g_n(t)\right]$, will have the form
\begin{align*}
U(t):=&\left[u_1,\ldots,u_{m_1},\ldots,u_{\left(\sum_{k=1}^{i-1} m_k\right)},\ldots,u_{\left(\sum_{k=1}^{i} m_k\right)},\ldots,u_{\left(\sum_{k=1}^d m_k\right)-1},u_{n}\right]\\
=&\begin{small}
\left[
  \begin{array}{ccc}
    \left[
      \begin{array}{cccc}
        v_{1,1} & v_{1,2} & \dots  & v_{1,m_1} \\
                & v_{1,1} & \ddots & v_{1,m_1-1} \\
                &         & \ddots & \vdots \\
                &         &        & v_{1,1} \\
      \end{array}
    \right]_{m_1\times m_1}&   &    \\
      & \ddots    &   \\
      &   &   \left[
      \begin{array}{cccc}
        v_{d,1} & v_{d,2} & \dots  & v_{d,m_d} \\
                & v_{d,1} & \ddots & v_{d,m_d-1} \\
                &         & \ddots & \vdots \\
                &         &        & v_{d,1} \\
      \end{array}
    \right]_{m_d\times m_d} \\
  \end{array}
\right]_{n\times n.}\end{small}
\end{align*}
For brevity, we focus on the $m_i$ elements contained in our $m_i\times 1$ solution vector $\nu$ with the understanding that it will be embedded into an $n\times 1$ vector of the form~\eqref{vec form u_i}, in the appropriate place, so that it only acts on the corresponding block matrix $K_i(t)$ of $K(t)$ in~\eqref{dyn eval prob for nu_i}.
There are $m_i$ linearly independent solutions of~\eqref{form of K} and they have the form
    \begin{align*}
    \nu_{i,1}(t)&:=\left[v_{i,m_i}(t),0,\ldots,0\right]_{m_i\times 1}^T\\
    \nu_{i,2}(t)&:=\left[v_{i,m_i-1}(t),v_{i,m_i}(t),0,\ldots,0\right]_{m_i\times 1}^T\\
    &\ \vdots\\
    \nu_{i,m_i-1}(t)&:=\left[v_{i,2}(t),v_{i,3}(t),\ldots,v_{i,m_i-1}(t),v_{i,m_i}(t),0\right]_{m_i\times 1}^T\\
    \nu_{i,m_i}(t)&:=\left[v_{i,1}(t),v_{i,2}(t),v_{i,3}(t),\ldots,v_{i,m_i-1}(t),v_{i,m_i}(t)\right]_{m_i\times 1}^T,
    \end{align*}
where the scalar functions $v_{i,j}$ are constructed by back solving based on the form of $\tilde{K}_i$ in \eqref{form of K}.  The associated dynamic equations are
    \begin{align}\label{assoc dyn eqs}
    v_{i,m_i}^\D(t)&=0\\
    v_{i,m_i-1}^\D(t)\notag&=\frac{1}{\l_i p}v_{i,m_i}(t)\\
    v_{i,m_i-2}^\D(t)\notag&=\frac{\m(t)-p}{2\l_i^2 p^2}v_{i,m_i}(t)+\frac{1}{\l_i p}v_{i,m_i-1}(t)\\
    v_{i,m_i-3}^\D(t)\notag&=\frac{(\m(t)-p)(\m(t)-2p)}{3!\l_i^3 p^3}v_{i,m_i}(t)+\frac{\m(t)-p}{2\l_i^2p^2}v_{i,m_i-1}(t)
        +\frac{1}{\l_i p}v_{i,m_i-2}(t)\\\notag
    &\ \vdots\\\notag
    v_{i,2}^\D(t)\notag&=\frac{(\m(t)-p)(\m(t)-2p)\dots(\m(t)-(m_i-3)p)}{(m_i-2)!\l_i^{m_i-2}p^{m_i-2}}v_{i,m_i}(t)\\\notag
    &+\frac{(\m(t)-p)(\m(t)-2p)\dots(\m(t)-(m_i-4)p)}{(m_i-3)!\l_i^{m_i-3}p^{m_i-3}}v_{i,m_i-1}(t)
        +\dots+\frac{\m(t)-p}{2\l_i^2 p^2}v_4(t)+\frac{1}{\l_i p}v_{i,3}(t)\\
    v_{i,1}^\D(t)\notag&=\frac{(\m(t)-p)(\m(t)-2p)\dots(\m(t)-(m_i-2)p)}{(m_i-1)!\l_i^{m_i-1}p^{m_i-1}}v_{i,m_i}(t)\\\notag
    &+\frac{(\m(t)-p)(\m(t)-2p)\dots(\m(t)-(m_i-3)p)}{(m_i-2)!\l_i^{m_i-2}p^{m_i-2}}v_{i,m_i-1}(t)
        +\dots+\frac{\m(t)-p}{2\l_i^2 p^2}v_{i,3}(t)+\frac{1}{\l_i p}v_{i,2}(t).\
    \end{align}
The solutions to~\eqref{assoc dyn eqs} are
    \begin{align*}
    v_{i,m_i}(t)&=1\\
    v_{i,m_i-1}(t)&=\int_{t_0}^t\frac{1}{\l_i p}v_{i,m_i}(\t)\D\t\\
    v_{i,m_i-2}(t)&=\int_{t_0}^t\frac{\m(\t)-p}{2\l_i^2 p^2}v_{i,m_i}(\t)\D\t
        +\int_{t_0}^t\frac{1}{\l_i p}v_{i,m_i-1}(\t)\D\t\\
    v_{i,m_i-3}(t)&=\int_{t_0}^t\frac{(\m(\t)-p)(\m(\t)-2p)}{3!\l_i^3p^3}v_{i,m_i}(\t)\D\t
        +\int_{t_0}^t\frac{\m(\t)-p}{2\l_i^2 p^2}v_{i,m_i-1}(\t)\D\t\\
        &+\int_{t_0}^t\frac{1}{\l_i p}v_{i,m_i-2}(\t)\D\t\\
    &\ \vdots\\
    v_{i,2}(t)&=\int_{t_0}^t\frac{(\m(\t)-p)(\m(\t)-2p)\dots(\m(\t)-(m_i-3)p)}{(m_i-2)!\l_i^{m_i-2}p^{m_i-2}}v_{i,m_i}(\t)\D\t\\
    &+\int_{t_0}^t\frac{(\m(\t)-p)(\m(\t)-2p)\dots(\m(\t)-(m_i-4)p)}{(m_i-3)!\l_i^{m_i-3}p^{m_i-3}}v_{i,m_i-1}(\t)\D\t\\
    &+\dots+\int_{t_0}^t\frac{\m(\t)-p}{2\l_i^2 p^2}v_{i,4}(\t)\D\t+\int_{t_0}^t\frac{1}{\l_i p}v_{i,3}(\t)\D\t\\
    v_{i,1}(t)&=\int_{t_0}^t\frac{(\m(t)-p)(\m(t)-2p)\dots(\m(t)-(m_i-2)p)}{(m_i-1)!\l_i^{m_i-1}p^{m_i-1}}v_{i,m_i}(t)\D\t\\
    &+\int_{t_0}^t\frac{(\m(\t)-p)(\m(\t)-2p)\dots(\m(\t)-(m_i-3)p)}{(m_i-2)!\l_i^{m_i-2}p^{m_i-2}}v_{i,m_i-1}(\t)\D\t\\
        &+\dots+\int_{t_0}^t\frac{\m(\t)-p}{2\l_i^2 p^2}v_{i,3}(t)\D\t+\int_{t_0}^t\frac{1}{\l_i p}v_{i,2}\D\t.
    \end{align*}


We can explicitly bound each $v_{i,j}$ which yields an explicit bound on each $u_{\ell_i+j}$, for all $1\leq j\leq m_i$ and for all $\ell_i$, with $1\leq i\leq d$.  Recall $\mu(t)\leq\m_{\max}\leq p$ for all $t\in\T$.  There exists constants $B_{i,j}$, $i=1,\dots,d$, and $j=1,\dots,m_i$, such that
    \begin{align*}
    |v_{m_i}(t)|&=1\leq B_{i,m_i}h_0(t,t_0)=B_{i,m_i},\\
    |v_{{m_i}-1}(t)|&\leq\int_{t_0}^t\frac{1}{\l_i p}v_{m_i}(\t)\D\t=\frac{h_1(t,t_0)}{\l_i p}\leq B_{i,{m_i}-1}h_1(t,t_0),\\
    |v_{{m_i}-2}(t)|&\leq\int_{t_0}^t\left|\frac{\m(\t)-p}{2\l_i^2 p^2}v_{m_i}(\t)\right|\D\t
    +\int_{t_0}^t\left|\frac{1}{\l_i p}v_{{m_i}-1}(\t)\right|\D\t\\
    &\leq\int_{t_0}^t\frac{p}{2\l_i^2 p^2}\D\t+\int_{t_0}^t\frac{h_1(\t,t_0)}{\l_i^2 p^2}\D\t\\
    &\leq \frac{h_1(t,t_0)}{2\l_i^2 p}+\frac{h_2(t,t_0)}{\l_i^2 p^2}\\
    &\leq B_{i,{m_i}-2}\sum_{j=1}^2h_j(t,t_0),\\
    |v_{{m_i}-3}(t)|&\leq\int_{t_0}^t\left|\frac{(\m(\t)-p)(\m(\t)-2p)}{3!\l_i^3 p^3}v_{m_i}(\t)\right|\D\t+\int_{t_0}^t\left|\frac{\m(\t)-p}{2\l_i^2 p^2}v_{{m_i}-1}(\t)\right|\D\t
    +\int_{t_0}^t\left|\frac{1}{\l_i p}v_{{m_i}-2}(\t)\right|\D\t\\
    &\leq\int_{t_0}^t\left|\frac{2p^2}{3!\l_i^3 p^3}\right|\D\t+\int_{t_0}^t\left|\frac{p}{2\l_i^2 p^2}\frac{h_1(\t,t_0)}{\l_i p}\right|\D\t
    +\int_{t_0}^t\left|\frac{1}{\l_i p}\left(\frac{h_1(\t,t_0)}{2\l_i^2 p}
    +\frac{h_2(\t,t_0)}{\l_i^2 p^2}\right)\right|\D\t\\
    &\leq B_{i,{m_i}-3}\sum_{j=1}^3h_j(t,t_0),\\
    |v_{{m_i}-4}(t)|&\leq\int_{t_0}^t\left|\frac{(\m(\t)-p)(\m(\t)-2p)(\m(\t)-3p)}{4!\l_i^4 p^4}v_{m_i}(\t)\right|\D\t+
    \int_{t_0}^t\left|\frac{(\m(\t)-p)(\m(\t)-2p)}{3!\l_i^3 p^3}v_{{m_i}-1}(\t)\right|\D\t\\
    &+\int_{t_0}^t\left|\frac{\m(\t)-p}{2\l_i^2 p^2}v_{{m_i}-2}(\t)\right|\D\t
    +\int_{t_0}^t\left|\frac{1}{\l_i p}v_{{m_i}-3}(\t)\right|\D\t\\
    &\leq\int_{t_0}^t\left|\frac{3!p^3}{4!\l_i^4 p^4}\right|\D\t+
    \int_{t_0}^t\left|\frac{2p^2}{3!\l_i^3 p^3}\frac{h_1(\t,t_0)}{\l_i p}\right|\D\t+\int_{t_0}^t\left|\frac{p}{2\l_i^2 p^2}\left(\frac{h_1(\t,t_0)}{2\l_i^2 p}
    +\frac{h_2(\t,t_0)}{\l_i^2 p^2}\right)\right|\D\t\\
    &+\int_{t_0}^t\left|\frac{1}{\l_i p}\left(\frac{h_1(\t,t_0)}{3\l_i^3 p}+2\frac{h_2(\t,t_0)}{2\l_i^3 p^2}+\frac{h_3(\t,t_0)}{\l_i^3 p^3}\right)\right|\D\t\\
    &\leq B_{i,{m_i}-4}\sum_{j=1}^4h_j(t,t_0),\\
    &\vdots\\
    |v_2|&\leq B_{i,2}\sum_{j=1}^{m_i-2}h_j(t,t_0),\\
    |v_1|&\leq B_{i,1}\sum_{j=1}^{m_i-1}h_j(t,t_0).
    \end{align*}
For each $1\leq i\leq d$, set $\beta_i:=\max_{j=1,\ldots,m_i}\{B_{i,j}\}.$
Then, for $1\leq i\leq d$ and $j=1,\ldots,m_i$, we have
$$
||u_{l_i+j}(t)||\leq \beta_i\sum_{k=0}^{m_i-1}h_k(t,t_0).
$$

Finally, recall that $w_i=Cu_i$ using $C$ from \eqref{doh} and define $D_i:=||C||\,\b_i$, for all $1\leq i\leq n$. Then we have the bound
$$
||w_i(t)||=||C u_i(t)||\leq ||C||\,\beta_i\sum_{k=0}^{m_i-1}h_k(t,t_0)=D_i\sum_{k=0}^{m_i-1}h_k(t,t_0),
$$
which proves the second claim in the theorem.
\end{proof}

To complete the setup for the main theorem, we state a necessary result from~\cite{Poochy} which characterizes the growth rates of generalized polynomials.

\begin{lemma}\label{lem-Poochy exp takes over polynomial}
Given a regressive scalar function $\l:\T\to\C$, $\l\in\textup{C}^1_\textup{rd}(\T,\C)$ on a time scale with bounded graininess and unbounded above with the property
    $$
    \exists T\in\T:0 < {-\inf}_{t\in [T,\infty)_{\T}} \left[\textup{Re}_\m (\l(t))\right],
    $$
it holds that
    $$
    \lim_{t\to\infty}h_k(t,t_0)e_\l(t,\t)=0,\quad\text{for }\t\in\T,\: k\in\N_0.
    $$
\end{lemma}

We are now in position to prove the main stability result.  The result will show that given a $p$-periodic time varying linear dynamic system~\eqref{Aprob}, the traditional eigenvalues of the associated time varying linear dynamic system~\eqref{Rprob} (via the Floquet decomposition of $\Phi_A$) completely determine the stability characteristics of the original system.

\begin{theorem}[{\bf Floquet Stability Theorem}]\label{Floquet_stab}
Given the $p$-periodic system~\eqref{Rprob} with eigenvalues $\{\g_i(t)\}_{i=1}^n$, we have the following properties of the solutions to the original $p$-periodic system~\eqref{Aprob}\textup{:}
\begin{enumerate}
\item[(i)] if Re$_\m\g_i(t)<0$ for all $i=1,\ldots,n$, then all the system~\eqref{Aprob} is exponentially stable;
\item[(ii)] if Re$_\m\g_i(t)=0$ for all $i=1,\ldots,n$, then in the case when only simple elementary divisors correspond to $\g_i(t)$, the system~\eqref{Aprob} is stable.  The system~\eqref{Aprob} is unstable, growing at rates of generalized polynomials of $t$ if among the elementary divisors there are multiple ones;
\item[(iii)] if Re$_\m\g_i(t)>0$ for some $i=1,\ldots,n$, then the system~\eqref{Aprob} is unstable.
\end{enumerate}
%
%
\end{theorem}

\begin{proof}
Let $\Phi_A(t,t_0)$ be the transition matrix for the system~\eqref{Aprob}.  Define the matrix $R(t)$ as in~\eqref{Rdef}, and note that by Theorem~\ref{Floquet2},~\eqref{Rprob} is the associated system through the Floquet decomposition of $\Phi_A$.  Given the traditional eigenvalues $\{\g_i(t)\}_{i=1}^n$ of $R(t)$, we define a set of dynamic eigenpairs of $R(t)$ by $\{\g_i(t),w_i(t)\}_{i=1}^n$ as in Lemma~\ref{Lem Existence of Dyn Epairs}.  By Theorem~\ref{thm-bd on dyn evecs}, the dynamic eigenvectors $w_i(t)$ are bounded by a finite sum of generalized polynomials as in~\eqref{dyn evec bound}.  Employing Theorem~\ref{mode vector stability theorem}, we can write the transition matrix of~\eqref{Rprob} as
    $$
    e_R(t,t_0)=\sum_{i=1}^n m_i(t) v_i^T(t_0),
    $$
for all $1\leq i \leq n$, where the vectors $m_i(t)$ and $v_i^T(t_0)$ are defined as in Theorem~\ref{mode vector stability theorem}.

Case (i): Applying Lemma~\ref{lem-Poochy exp takes over polynomial}, for each $1\leq i \leq n$, we have
    $$
    \lim_{t\to\infty} ||m_i(t)||\leq D_i\sum_{k=0}^{m_i-1}h_k(t,t_0)|e_{\g_i}(t,t_0)|\leq C_\e e_{\textup{Re}_\m[{\g_i}]\oplus\e}(t,t_0)
    \leq C_\e e_{-\l}(t,t_0)=0,
    $$
where $\e>0$ is chosen so that $\g_i(t)\oplus\e\in\H_t$, for all $t\in\T$, $C_\e:=\max_{t\geq t_0} D_i\sum_{k=0}^{m_i-1}h_k(t,t_0)e_{\ominus\e}(t,t_0)$, and $\l>0$ with Re$_\m[\g_i(t)]\oplus\e<-\l<0$, for all $i=1,\ldots,n$ and $t\in\T$.
By Theorem~\ref{mode vector stability theorem} (iii), the system~\eqref{Rprob} is exponentially stable.

The solution to~\eqref{Aprob} and the solution to~\eqref{Rprob} are related by a Lyapunov transformation, namely $L(t):=\Phi_A(t,t_0)e_R^{-1}(t,t_0)$.  By Theorem~\ref{relationship_thm} and Corollary~\ref{cor stability thru Lyap trans}, the solution to the system~\eqref{Aprob} is exponentially stable if and only if the solution to the system~\eqref{Rprob} is exponentially stable.

Case (ii): Suppose Re$_\m[\g_i(t)]=0$ for all $i=1,\ldots,n$ with only simple elementary divisors corresponding to each $\g_i(t)$.  Then each Jordan block corresponding to each $\g_i(t)$ is $1\times 1$, which implies that the mode vectors have the form
    $$
    m_i(t)=\b_i e_{\g_i}(t,t_0).
    $$
Thus,
    $$
    \lim_{t\to\infty} ||m_i(t)||\leq \lim_{t\to\infty} \b_i e_{\textup{Re}_\m[{\g_i}]}(t,t_0)=\b_i<\infty.
    $$
By Theorem~\ref{mode vector stability theorem} (ii), the system~\eqref{Rprob} is stable.  Thus, since~\eqref{Aprob} is related to~\eqref{Rprob} by a Lyapunov transformation,~\eqref{Aprob} is stable.

Now suppose that for some $i$ the elementary divisor corresponding to $\g_i(t)$ has multiplicity $n>1$.  Then, the elements in the corresponding mode vector is made up of generalized polynomials.  Thus, since there is no decaying exponential function to control the growth of these polynomials, solutions grow at the rate of the polynomials.  Therefore,~\eqref{Rprob} is unstable, implying by the Lyapunov transformation that~\eqref{Aprob} is also unstable.

Case (iii): Suppose that for some $1\leq i\leq n$, we have Re$_\m[\g_i(t)]>0$.  This implies that
    $$
    \lim_{t\to\infty}||e_R(t,t_0)||=\infty,
    $$
thus implying by the Lyapunov transformation
    $$
    \lim_{t\to\infty}||\Phi_A(t,t_0)||=\infty.
    $$
\end{proof}

The consequences of Theorem~\ref{Floquet_stab} should be
underscored at this point:  from the construction of the
corresponding system matrix $R(t)$, we can deduce the stability
of~\eqref{Rprob} solely on the placement of the traditional
eigenvalues in the complex plane.  In general, this is not true
of time varying systems.  However, we have shown that for this
unified Floquet theory, the associated (but in general time
varying) matrix $R(t)$ yields stability characteristics through
simple pole placement in the complex plane---exactly like
constant systems and Jordan reducible time varying systems on
$\R$ or $h\Z$.  Moreover, the matrix $R(t)$ is in fact constant
if and only if the domain of the underlying dynamical system has
constant graininess, e.g. $\R$ or $h\Z$. The upshot of course is
that this Floquet theory extends to any domain that is a
$p$-periodic closed subset of $\R$, including those that have
nonconstant discrete points as well as mixed continuous and
discrete intervals.

A key observation about this unified Floquet theory is that the
associated system matrix $R(t)$ yields the stability
characteristics of the original system by placement of the
traditional eigenvalues alone---although it accomplishes this
via an associated time varying systems rather than with a time
invariant system.

\begin{remark}
We can think of an eigenvalue
$\g(t)=\lim_{s\searrow\m(t)}\frac{\l^{\frac{s}{p}}-1}{s}$ of
$R(t)$ as constant with respect to its ``placement'' in the
instantaneous stability region corresponding to the element $t$
in the  $p$-periodic time scale $\T$.  In other words, the
exponential bound is independent of $\mu(t)$:
    \begin{align*}
    e_{\g}(t,t_0)&=\exp\left[\int_{t_0}^t\lim_{s\searrow\m(\t)}\Log\left(1+s\g(\t)\right)\D\t\right]\\
    &=\exp\left[\int_{t_0}^t\lim_{s\searrow\m(\t)}\Log\left(1+s\frac{\l^{\frac{s}{p}}-1}{s}\right)\D\t\right]\\
    &=\l^{\frac{t-t_0}{p}}.
    \end{align*}
%
\end{remark}

\section{Conclusions}
A unified and extended Floquet theory has been developed, which
includes a canonical Floquet decomposition in terms of a generalized
matrix exponential function.  The notion of time varying Floquet
exponents and their relationship to constant Floquet
multipliers has been introduced.  We presented a closed-form
answer to the open problem of finding a solution matrix $R(t)$ to
the equation $e_R(t,\t)=M$ on an arbitrary time scale $\T$, where $R(t)$ and $M$ are $n\times
n$ matrices, and $M$ is constant and nonsingular. The development of
this Floquet theory on both homogeneous and nonhomogeneous hybrid
dynamical systems was completed by developing Lyapunov
transformations and a stability analysis of the original system in
terms of a corresponding (but not necessarily autonomous) system
that, by its construction, yields stability
characteristics by traditional eigenvalue placement in the complex plane.


\section{Appendix: Real Powers of a Matrix}\label{Appendix-Real Pwr of a Mtrx}

We introduce the technical preliminaries required to define the
real power of a matrix, which is of paramount importance in the
derivation of the unified Floquet theory. The following can be
found in \cite{ho1, ho2}.

\begin{definition}\label{ho's projection matrix P(A) definition}
Given an $n\times n$ invertible matrix $M$ with elementary divisors $\{(\l-\l_i)^{m_i}\}_{i=1}^k$, we let the
characteristic polynomial be the function $p(\l)$.  We define
the polynomial $a_i(\l)$ implicitly via the expansion
    $$
    \frac{1}{p(\l)}=\sum_{i=1}^k\frac{a_i(\l)}{(\l-\l_i)^{m_i}}.
    $$
The polynomial $b_i(\l)$ is defined by omitting the factor
$(\l-\l_i)^{m_i}$ from the characteristic polynomial $p(\l)$,
that is
    $$
    b_i(\l):=\prod_{j\neq i}(\l-\l_j)^{m_j}.
    $$
Then the {\it $i$th projection matrix} is defined as
    \begin{equation}\label{ho's projection matrix P(A) equation}
    P_i(\l):=a_i(\l)b_i(\l).
    \end{equation}
\end{definition}

A property of the set of the projection matrices is the
following.  By definition,
    \begin{equation}\label{scalar identity prop of projection
    matrix}
    1 =\sum_{j=1}^k\frac{a_j(\l)}{(\l-\l_j)^{m_j}}p(\l)=\sum_{j=1}^ka_j(\l)b_j(\l)=\sum_{j=1}^kP_j(\l).
    \end{equation}
When the scalar value $\l$ is replaced by the $n\times n$
nonsingular matrix $M$, \eqref{scalar
identity prop of projection matrix} becomes
    \begin{equation}\label{matrix identity prop of projection
    matrix} I =\sum_{j=1}^kP_j(M).
    \end{equation}

We now state and prove three propositions concerning the
projection matrices.

\begin{prop}\label{proj matrix proposition 1}
Given an $n\times n$ nonsingular matrix $M$ as in Definition~\ref{ho's projection matrix P(A) definition}, we have
$
P_i(M)(M-\l_iI)^{m_i}=0.
$
\end{prop}

\begin{proof}
Observe
$$
P_i(M)(M-\l_iI)^{m_i}=a_i(M)b_i(M)(M-\l_iI)^{m_i}=a_i(M)\prod_{j=1}^k(M-\l_jI)^{m_j}=a_i(M)p(M)=0.
$$
\end{proof}

\begin{prop}\label{proj matrix proposition 2-proj matrices are orthog}
The set $\{P_i(M)\}_{i=1}^k$ is orthogonal.  That is,
$P_i(M)P_j(M)=\delta_{ij}P_i(M)$.
\end{prop}

\begin{proof}
Suppose $i\neq j$. Then
    \begin{align*}
    P_i(M)P_j(M)&=P_i(M)a_j(M)b_j(M)\\
    &=P_i(M)a_j(M)\prod_{k\neq j}(M-\l_kI)^{m_k}\\
    &=P_i(M)a_j(M)(M-\l_iI)^{m_i}\prod_{k\neq i,j}(M-\l_kI)^{m_k}\\\\
    &=P_i(M)(M-\l_iI)^{m_i}a_j(M)\prod_{k\neq i,j}(M-\l_kI)^{m_k}\\
    &=0.
    \end{align*}
To conclude the proof, using the property~\eqref{matrix identity
prop of projection matrix} and multiplying on each side by the
$i$th projection matrix, we have
    $$
    P_i(M)=\sum_{j=1}^kP_i(M)P_j(M)=P_i(M)P_i(M).
    $$
\end{proof}

\begin{prop}\label{prop-proj matrix invariance with JCF}
Given an $n\times n$ nonsingular matrix $M$ as in Definition~\ref{ho's projection matrix P(A) definition} with corresponding Jordan canonical form $J=CMC^{-1}$, for some nonsingular matrix $C$, we have
$$
CP_i(M)C^{-1}=P_i(J),
$$
for all $1\leq i\leq k$.
\end{prop}
\begin{proof}
For any $1\leq i\leq k$, by definition of the projection matrix~\eqref{ho's projection matrix P(A) equation}, we see that $P_i(M)$ is the product of two polynomials of the matrix $M$.  Further, by definition of $a_i$ and $b_i$, it follows that $Ca_i(M)C^{-1}=a_i(J)$ and $Cb_i(M)C^{-1}=b_i(J)$.  Thus, $$CP_i(M)C^{-1}=Ca_i(M)b_i(M)C^{-1}=Ca_i(M)C^{-1}Cb_i(M)C^{-1}=a_i(J)b_i(J)=P_i(J).$$
\end{proof}

We are now in the position to state the definition of the principal value of the real
power of a nonsingular $n\times n$ matrix.
\begin{definition}\label{matrix M to a real power definition-APPENDIX DEFINITION}
Given an $n\times n$ nonsingular matrix $M$ as in Definition~\ref{ho's projection matrix P(A) definition} and any $r\in\R$, we
define the real power of the matrix $M$ by
    $$
    M^r:=\sum_{i=1}^kP_i(M)\l_i^r\left[\sum_{j=0}^{m_i-1}\frac{\Gamma(r+1)}{j!\;\Gamma(r-j+1)}\left(\frac{M-\l_iI}{\l_i}\right)^j\right].
    $$
\end{definition}

The next theorem gives a relationship between the real power of a matrix and the real power of the corresponding Jordan canonical form.
\begin{theorem}\label{Thm-real power of M relationship with real power of JCF}
Given an $n\times n$ nonsingular matrix $M$  as in Definition~\ref{ho's projection matrix P(A) definition}with corresponding Jordan canonical form $J=CMC^{-1}$, for an appropriately defined nonsingular matrix $C$, and any $r\in\R$, we have
$$
CM^rC^{-1}=J^r.
$$
\end{theorem}
\begin{proof}
The result follows from a direct application of Proposition~\ref{prop-proj matrix invariance with JCF} and Definition~\ref{matrix M to a real power definition-APPENDIX DEFINITION}.  Observe,
\begin{align*}
CM^rC^{-1}&=C\sum_{i=1}^kP_i(M)\l_i^r\left[\sum_{j=0}^{m_i-1}
\frac{\Gamma(r+1)}{j!\;\Gamma(r-j+1)}\left(\frac{M-\l_iI}{\l_i}\right)^j\right]C^{-1}\\
&=\sum_{i=1}^kCP_i(M)C^{-1}C\l_i^r\left[\sum_{j=0}^{m_i-1}
\frac{\Gamma(r+1)}{j!\;\Gamma(r-j+1)}\left(\frac{M-\l_iI}{\l_i}\right)^j\right]C^{-1}\\
&=\sum_{i=1}^kP_i(J)\l_i^r\left[\sum_{j=0}^{m_i-1}
\frac{\Gamma(r+1)}{j!\;\Gamma(r-j+1)}\left(\frac{J-\l_iI}{\l_i}\right)^j\right]\\
&=J^r.
\end{align*}
\end{proof}

\begin{remark}
Theorem~\ref{Thm-real power of M relationship with real power of JCF} is a generalization from the familiar matrix property that $CM^kC^{-1}=J^k$, for all $k\in\N_0$.
\end{remark}

\begin{remark}\label{nilpotency remark}
Given any nonsingular $n\times n$ matrix $M$ as in Definition~\ref{ho's projection matrix P(A) definition}, we have $P_i(M)(M-\l_iI)^{m_i}=0$, for $i=1,\ldots,k$ and $m_1+\cdots+m_k=n$.  However, since
$P_i(M)(M-\l_iI)$ is nilpotent by Proposition~\ref{proj matrix
proposition 1} and Proposition~\ref{proj matrix proposition
2-proj matrices are orthog}, there exists an integer $n_i \leq
m_i$ such that $P_i(M)(M-\l_iI)^{n_i-1}\neq 0$, but
$P_i(M)(M-\l_iI)^{n_i}=0$. Then $P_i(M)(M-\l_iI)^{n_i}$ is a
minimal polynomial of $M$.
\end{remark}

\begin{theorem}\label{ho's eigenvalue thm}
A nonzero column vector $v_r$ of $P_i(M)(M-\l_iI)^{n_i-r}$ is a
generalized eigenvector of rank $r$, for $1\leq r\leq n_i$,
associated with the eigenvalue $\l_i$.  A regular eigenvector
has rank $r=1$.
\end{theorem}
\begin{proof}
For $r=1$, we have
\begin{align*}
MP_i(M)(M-\l_iI)^{n_i-1}&=\sum_{j=1}^kP_j(M)\left[\l_jI+(M-\l_jI)\right]P_i(M)(M-\l_iI)^{n_i-1}\\
&=\l_iP_i(M)(M-\l_iI)^{n_i-1},
\end{align*}
where we use the result from Proposition~\ref{proj matrix proposition 2-proj matrices are orthog} which states that $P_i(M)P_j(M)=\delta_{ij}P_i(M)$, for all $1\leq i,j\leq k$.
Thus, since $P_i(M)(M-\l_iI)^{n_i-1}$ has rank one, we choose a
nonzero $n\times 1$ column vector $v_1$ from this matrix and
obtain
$$
Mv_1=\l_iv_1.
$$
When $r>1$, note that
$$
P_i(M)(M-\l_iI)^{n_i-r}=(M-\l_iI)P_i(M)(M-\l_iI)^{n_i-r-1,}
$$
which implies
$$
v_r=(M-\l_iI)v_{r+1}.
$$
\end{proof}

The following theorem will show that for any invertible $n\times
n$ matrix $M$, with $d$ eigenpairs $\{\l_i,v_i\}_{i=1}^d$, and
any $r\in\R$, the matrix $M^r$ has the $d$ eigenpairs
$\{\l_i^r,v_i\}_{i=1}^d$.

\begin{theorem}\label{matrix M to a real power theorem}
Given an invertible $n\times n$ matrix $M$ as in Definition~\ref{ho's projection matrix P(A) definition}, with $d$ eigenpairs
$\{\l_i,v_i\}_{i=1}^d$, and any $r>0$, for each eigenpair
$\{\l_i,v_i\}$ of $M$, we have $M^rv_i=\l_i^rv_i$.
\end{theorem}

\begin{proof}
Let $1\leq s\leq d$ and $\l_s$ be an eigenvalue of the matrix
$M$. Then, noting that $(M-\l_sI)^{n_s}$ is a minimal polynomial of $M$, using Proposition~\ref{proj matrix proposition 1} and Remark~\ref{nilpotency remark} we have
    \begin{align*}
    M^rP_s(M)&(M-\l_sI)^{n_s-1}\\
    &=\left[\sum_{i=1}^dP_i(M)\l_i^r\left[\sum_{j=0}^{m_i-1}
    \frac{\Gamma(r+1)}{j!\;\Gamma(r-j+1)}\left(\frac{M-\l_iI}{\l_i}\right)^j\right]\right]P_s(M)(M-\l_sI)^{n_s-1}\\
    &=P_s(M)\l_s^r\sum_{j=0}^{m_i-1}
    \frac{\Gamma(r+1)}{j!\;\Gamma(r-j+1)}\l_i^{-j}\left(M-\l_iI\right)^{j+n_s-1}\\
    &=\l_s^rP_s(M)(M-\l_sI)^{n_s-1},
    \end{align*}
for each nonzero column vector $v_s$ of
$P_s(M)(M-\l_sI)^{n_s-1}$.  The matrix $P_s(M)(M-\l_sI)^{n_s-1}$
has rank 1 and we can take any nonzero column vector $v_s$ from
this matrix as an eigenvector corresponding the eigenvalue
$\l_s$ of the matrix $M$.  This vector is a linear combination
of all of the eigenvectors associated with the eigenvalue
$\l_s$. Thus, $M^rv_s=\l_s^rv_s$.
\end{proof}

\end{document}